\documentclass[reqno,11pt]{amsart}
\usepackage{amsfonts,amsmath,amssymb}
\usepackage{graphicx}

\usepackage{hyperref}

\usepackage[margin=1.23in]{geometry}
\numberwithin{equation}{section}

\def\R{{\mathbb R}}

\def\e{{\varepsilon}}

\def\eps{{\varepsilon}}
\def\varep{{\varepsilon}}

\def\+R{+_{_{ \!\! \R}}}

\def\div{{\,\mbox{div}\,}}
\def\curl{{\,\mbox{curl}\,}}

\def\bar{\overline}

\DeclareMathAlphabet{\mathpzc}{OT1}{pzc}{m}{it}

\numberwithin{equation}{section}

\theoremstyle{plain}

\newtheorem{theorem}{Theorem}[section]

\newtheorem{proposition}[theorem]{Proposition}

\begin{document}

\title[Global regularity for water waves]{Recent advances on the global regularity \\ for irrotational water waves}

\author{A. D. Ionescu}
\address{Princeton University}
\email{aionescu@math.princeton.edu}

\author{F. Pusateri}
\address{Princeton University}
\email{fabiop@math.princeton.edu}

\thanks{The first author was supported in part by NSF grant DMS-1600028 and by NSF-FRG grant DMS-1463753.
The second author was supported in part by NSF Grant DMS-1265875.}

\begin{abstract}
We review recent progress on the long-time regularity of solutions of the Cauchy problem for the water waves equations, in two and three dimensions.

We begin by introducing the free boundary Euler equations and discussing the local existence of solutions using the paradifferential approach,
as in \cite{AlMet1,ABZ1,ABZ2}.
We then describe in a unified framework, using the Eulerian formulation,
global existence results for three dimensional and two dimensional gravity waves, see \cite{GMS2,Wu3DWW,WuAG,IoPu2,ADa,ADb,HIT,IT,Wa1}, and our joint result 
with Deng and Pausader \cite{DIPP} on global regularity for the 3D gravity-capillary model.

We conclude this review with a short discussion about the formation of singularities,
and give a few additional references to other interesting topics in the theory.

\end{abstract}

\maketitle

\begin{quote}
%\footnotesize
\setcounter{tocdepth}{1}
\tableofcontents
\end{quote}

\section{Introduction}

%\subsection{The water waves problem}\label{secWW0}
The study of the motion of water waves, such as those on the surface of the ocean, is a classical question,
and one of the main problems in fluid dynamics.
The origins of water waves theory can be traced back
at least to the work of Laplace and Lagrange, Cauchy \cite{CauchyMemoirs} and Poisson, %\cite{PoissonMemoirs},
and then Russel, Green and Airy, among others, see \cite{Craik}.
%The theory of water waves, and fluids more in general, has historically impacted greatly mathematical and analytical developments.
% Cauchy -> Complex function theory
% Cauchy + Poisson -> Fourier analysis
Classical studies include those by Stokes \cite{Stokes}, Levi-Civita \cite{LeviC}
and Struik \cite{Struik} on progressing waves, %in infinite depth / finite depth ,
the instability analysis of Taylor \cite{Taylor},
the works on solitary waves by Friedrichs and Hyers \cite{FriedHy}, %Beale \cite{Beale} and Amick and Toland \cite{AmTol1},
and on steady waves by Gerber \cite{Gerber}.

The main questions one can ask about water waves are the typical ones for any physical evolution problem:
the existence of solutions of the initial value problem, their regularity,
the possible formation of various singularities in the flow,
the existence of special solutions (such as solitary waves, standing waves, periodic/quasiperiodic waves) and their stability,
and the long-time existence and asymptotic behavior of the flow.
There is a vast body of literature dedicated to all of these aspects.

The main focus of this article is to review the local and global existence theory for the initial value problem
associated to the water waves equations, and give an overview of the recent progress in this area. We will refer the reader to various books, research papers and surveys for others aspects of the theory.

We will concentrate on the motion of an inviscid and irrotational %water waves problem
$2$ or $3$ dimensional fluid occupying a region of infinite depth and infinite extent below the graph of a function.
These are models for the motion of waves on the surface of the deep ocean,
where the two dimensional case corresponds to waves whose motion is assumed to be constant in one direction on the interface.
We will consider both $2$ and $3$ dimensional dynamics under the influence of the gravitational force
and/or surface tension effects acting on particles at the interface.
Our main goal is to present, in a unified framework, several results about the global existence of solutions which are initially small,
that is, sufficiently close to a flat and still interface
in a suitable sense.

\subsection{Structure of the paper}
%Before describing in more details these results, %concerning the global existence of solutions for the water waves problem,
In section \ref{secWW} we introduce the free boundary Euler equations in the standard Eulerian formulation
%\ref{secWW} below,
and the Zakharov-Craig-Schanz-Sulem Hamiltonian formulation for irrotational flows.
In section \ref{secLoc} we discuss the short time existence of solutions following the paradifferential
approach of \cite{AlMet1,ABZ1,ABZ2}.
%The remaining sections will then be dedicated to describing global existence results in both $3$ and $2$ dimensions.
Section \ref{global} is dedicated to global existence results.
We discuss three different problems, in increasing order of difficulty: the 3D gravity water waves, the 2D gravity water waves, and the 3D gravity-capillary water waves.
Section \ref{secOther} contains a brief discussion about the formation of singularities,
and few additional references to other interesting topics in the theory.

%In this introduction we will first describe the equations of motion and give a %statement of our main result, Theorem \ref{MainTheo}.
%Our main results are stated in \ref{MainResult}.
%We will then provide some history and background concerning the local wellposedness of %the problem and related questions.
%We will also discuss in some generality the question of global regularity, some %previous works on this problem, and some of our main ideas.
%Some of the main ideas behind our results are discussed in fairly general terms during %this introduction.
%The strategy of our proof is then described in more technical terms in section \ref{MPO}.
%A discussion of the key ideas involved in the various steps of the proof is given in %\ref{secideas}.

\subsection{Free boundary Euler equations}\label{secWW}
The evolution of an inviscid perfect fluid that occupies a domain $\Omega_t \subset \R^n$, for $n \geq 2$, at time $t \in \R$,
is described by the free boundary incompressible Euler equations.
We let $v$ and $p$ denote respectively the velocity and the pressure of the fluid, at time $t$ and position $x \in \Omega_t$,
and assume that the fluid has constant density equal to $1$.
If the fluid evolves in a gravitational field, the equations of motion are
\begin{subequations}
\label{Euler}
\begin{equation}\label{E}
\partial_t v + v \cdot \nabla v = - \nabla p - g e_n, \qquad \nabla \cdot v = 0, \qquad x \in \Omega_t,
\end{equation}
where $e_n$ is the $n$-th standard unit vector of $\R^n$ and $g$ is the gravitational constant.
The first equation in \eqref{E} is the conservation of momentum equation, while the second is the incompressibility condition.
When gravitational effects are neglected one sets $g=0$ in \eqref{E}.

The boundary of the fluid evolves with time and is part of the unknowns in the problem.
In  particular, the free surface $S_t := \partial \Omega_t$
moves with the normal component of the velocity according to the kinematic boundary condition
 \begin{equation}\label{BC1}
\partial_t + v \cdot \nabla  \qquad \mbox{is tangent to} \quad {\bigcup}_t S_t \subset \R^{n+1}_{t,x}.
\end{equation}
%%% USE t first then x.
%To complete the system one needs to add one more boundary condition for the pressure
The atmospheric pressure outside the fluid domain is assumed to be constant, and set to zero for convenience.
On the interface the pressure is given by
\begin{equation}
\label{BC2}
p (t,x) = \sigma \kappa(t,x),  \qquad x \in S_t,
\end{equation}
\end{subequations}
where $\kappa$ is the mean-curvature of $S_t$ and $\sigma \geq 0$ is the surface tension coefficient.
At liquid-air interfaces, the surface tension force results from the greater attraction
of water molecules to each other than to the molecules in the air.
%It is interesting to notice that water has the highest surface tension among all common liquids except mercury.\footnote{See
%\cite[sec. 2.4.1]{ConstBook} for a discussion about the nature of surface tension.}
%Surface tension effects are only relevant at short wavelength, tipically for waves that do not exceed a few centimeters.

%\begin{figure}[h!]
%%{\includegraphics[width=0.5\textwidth,height=0.2\textheight]{?.jpeg}} \\
%\tofill
%\caption{\small Put figure of general domain ? (to create)}
%  \label{fig:WW0}
%\end{figure}

One can consider the free boundary Euler equations above in various types of domains $\Omega_t$ (bounded, periodic, unbounded)
and study flows with different characteristics (rotational or irrotational, with gravity and/or surface tension),
or even more complicated scenarios where the moving interface separates two fluids.

There are several difficulties in treating the system \eqref{E}-\eqref{BC1}-\eqref{BC2}
which are due to the quasilinear nature of the equations, i.e., the highest derivatives appear nonlinearly,
and, above all, to the free moving boundary and its interaction with the fluid.
As we will discuss below, the system \eqref{E}-\eqref{BC1}-\eqref{BC2} has a %degenerate
``hyperbolic'' structure, which can only be captured thanks to great insights into the nature of the equations,
as was done for example in \cite{Wu1,Wu2,CL,Lindblad,Lannes,CS2,ShZ1,ABZ1}.
This structure leads to a priori control and local existence of solutions for sufficiently
regular initial data in the case of non-self intersecting interfaces, provided that
\begin{align}
 \label{RT}
-\nabla_N p(t,x) > 0, \qquad x \in S_t,
\end{align}
where $N$ is the outer unit normal, which is the so-called {\it Rayleigh-Taylor} sign condition.
In general, when \eqref{RT} is violated instabilities might occur \cite{Ebin1,BHL,BarLan}.

We will discuss in more details these local regularity issues in Section \ref{secLoc} below by restricting our attention to the case of irrotational flows,
and following the paradifferential approach of Alazard--M\'etivier \cite{AlMet1} %\cite{Alinhac}
and Alazard--Burq--Zuily \cite{ABZ1,ABZ2}.
We choose to present this approach among the various possible ones, since it is well-suited as a starting point
for the discussion of the long-time regularity results in Section \ref{global}. %\ref{sec3d}, \ref{sec2d} and \ref{sec3dfull}.

\medskip
\subsection{The Hamiltonian formulation} %Zakharov-Craig-Schanz-Sulem formulation}
\label{secWW2}
In the case of irrotational flows, that is when
\begin{align}
 \label{irr}
\rm{\curl} v = 0,
\end{align}
one can reduce \eqref{Euler} to a system of two equations on the boundary.
Indeed, assume also that $\Omega_t \subset \R^n$ is the region below the graph of a function
$h : I_t \times \R^{n-1}_x  \rightarrow \R$,
that is
\begin{align}
\label{Omegat}
\Omega_t = \{ (x,y) \in \R^{d} \times \R \, : y \leq h(t,x) \} \quad \mbox{and} \quad S_t = \{ (x,y) : y = h(t,x) \}, \qquad d:=n-1.
\end{align}
Let $\Phi$ denote the velocity potential,
\begin{align}
\label{Phi}
\nabla_{x,y} \Phi(t,x,y) = v (t,x,y), \qquad \Delta_{x,y} \Phi(t,x,y) = 0, \qquad (x,y) \in \Omega_t,
\end{align}
and let
\begin{align}
\label{phi}
\phi(t,x) := \Phi (t, x, h(t,x))
\end{align}
denote the restriction of $\Phi$ to the boundary $S_t$.

\begin{figure}[h!]
{\includegraphics[width=0.6\textwidth,height=0.25\textheight]{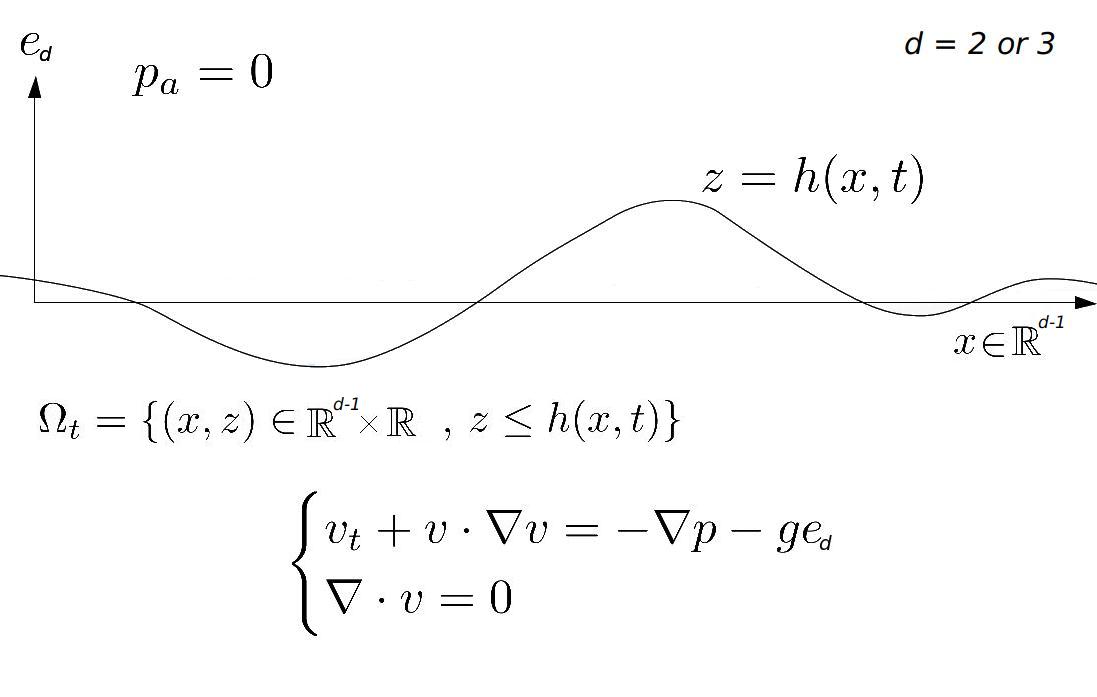}}\\
\caption{\small The water waves problem in Eulerian coordinates.}
  \label{fig:WW}
\end{figure}

Then, the equations of motion can be reduced to the following system for the unknowns
$h, \phi : I_t \times \R^{n-1}_x \rightarrow \R$:
\begin{equation}
\label{WWE}
\left\{
\begin{array}{l}
\partial_t h = G(h) \phi,
\\
\partial_t \phi = -g h  + \sigma \div \Big[ \dfrac{\nabla h}{ (1+|\nabla h|^2)^{1/2} } \Big]
  - \dfrac{1}{2} {|\nabla \phi|}^2 + \dfrac{{\left( G(h)\phi + \nabla h \cdot \nabla \phi \right)}^2}{2(1+{|\nabla h|}^2)}.
\end{array}
\right.
\end{equation}
Here
\begin{equation}
\label{defG0}
G(h) := \sqrt{1+{|\nabla h|}^2} \mathcal{N}(h), \qquad   \mathcal{N}(h) \phi := \nabla_n \Phi,
\end{equation}
where $n$ is the outward unit normal to $S_t$,
and $\mathcal{N}(h)$ denotes the Dirichlet--Neumann map associated to the domain $\Omega_t$.
Roughly speaking, one can think of $G(h)$ as a first order, non-local, linear operator that depends nonlinearly on the domain.
%and on a finite number of derivatives of $h$.
%For smooth and small $h$, a good approximation of it is given by $G(h)\phi = |\nabla| \phi$.
We refer to  \cite[chap. 11]{SulemBook} or the book of Lannes \cite{LannesBook} for the derivation of \eqref{WWE},
which is the so-called Zakharov-Craig-Schanz-Sulem \cite{CSS} formulation. %% Ref to Zakharov

The system \eqref{WWE} admits the conserved energy
\begin{equation}\label{CPWHam}
\begin{split}
\mathcal{H}(h,\phi) &:= \frac{1}{2} \int_{\R^{n-1}} G(h)\phi \cdot \phi \, dx + \frac{g}{2} \int_{\R^{n-1}} h^2 \,dx
  + \sigma\int_{\R^{n-1}} \frac{{|\nabla h|}^2}{1 + \sqrt{1+|\nabla h|^2} } \, dx
\end{split}
\end{equation}
which is the sum of the kinetic energy corresponding to the $L^2$ norm of the velocity field
and the potential energy due to gravity and surface tension.
It was first observed by Zakharov \cite{Zak0} that \eqref{WWE} is the Hamiltonian\footnote{Recently, Craig \cite{CraigHamWW}
has shown that \eqref{Euler}-\eqref{irr} can be formulated as an Hamiltonian system for general smooth domains.}
flow associated to \eqref{CPWHam}. For sufficiently small and smooth solutions one has
\begin{equation}\label{CPWHam'}
\begin{split}
\mathcal{H}(h,\phi) & \approx {\big\| |\nabla |^{1/2} \phi \big\|}_{L^2}^2 + {\big\| (g-\sigma\Delta)^{1/2}h \big\|}_{L^2}^2.
\end{split}
\end{equation}
%where $L^2$ denotes the standard space of square integrable functions.
%and \eqref{WWE} is usually referred to as the Craig-Schanz-Sulem-Zakharov formulation \cite{CSS}.

%This system describes the evolution of an incompressible perfect fluid of infinite depth and infinite extent,
%with a free moving surface which is the graph of a smooth function,
%and a pressure boundary condition given by the Young-Laplace equation.

The formal linearization of \eqref{WWE}-\eqref{defG0} around a flat and still interface is
\begin{align}
\label{WWElin}
\partial_t h = |\nabla| \phi, \qquad \partial_t \phi = -g h  + \sigma \Delta h.
\end{align}
By defining the linear dispersion relation
\begin{align}
\label{disprel}
\Lambda_{g,\sigma} := \sqrt{g|\nabla|+\sigma|\nabla|^3},
\end{align}
the identitites \eqref{WWElin} can be written as a single equation for a complex-valued unknown,
\begin{align}
\label{WWElin2}
\partial_t u + i \Lambda_{g,\sigma}u = 0, \qquad  u := \sqrt{g+\sigma|\nabla|^2}h + i|\nabla|^{1/2}\phi.
\end{align}

One generally refers to \eqref{WWE} as the {\bf gravity water waves} system when $g>0$ and $\sigma=0$,
as the {\bf capillary water waves} system when $g=0$ and $\sigma>0$,
and as the {\bf gravity-capillary water waves} system when $g>0$ and $\sigma>0$.

We remark that the presence of the various forces (gravity and/or surface tension)
does have an impact on the existence theory of solutions.
In the local existence theory this impact is mostly quantitative and %Better phrasing
the techniques developed for a specific scenario are likely to be adaptable to others.
On the other hand, when considering the long-time existence of solutions, the presence of gravity and/or surface tension has a major impact on the evolution.
This can be seen already at the level of the linearized equations \eqref{WWElin2},
and is even more apparent when looking at three waves resonant interactions, (quadratic time-resonances),
and at fully coherent interactions (space-time resonances).
We explain these concepts in Section \ref{sec3d}.

\subsection{Other approaches and formulations}\label{otherappr}

There are of course other possible descriptions of the equations.
In the series of works \cite{Wu1}--\cite{Wu3DWW}, see also the survey \cite{Wusurvey},
Wu uses a combination of Lagrangian coordinates and tools from complex analysis
such as the Riemann mapping Theorem and the theory of holomorphic functions (Clifford analysis in $3$d).
%complex analysis very effective in the two dimensional setting.
%Also \cite{IoPu4} for Dirichlet--Neumann in $1$d
Lagrangian coordinates and variations have been used in \cite{Lindblad,CL,CS2},
and complex coordinates in the works of Nalimov \cite{Nalimov},
in Zakharov et al. \cite{DKSZ} %Kuznetsov--Spector--Zakharov \cite{KSZ}, not this one?
in various theoretical and numerical works, see for example \cite{ZakDya2,DyaLuKo} and references therein, and in the series of papers \cite{HIT}-\cite{ITG}.
%by Hunter--Ifrim--Tataru \cite{HIT}, Ifrim--Tataru \cite{IT2}, and Harrop-Griffiths--Ifrim--Tataru \cite{BIT}. %Can skip this last here?
%There are various different approaches to the formulation of the water waves system.

%Local theory works well anyways, some approaches are more natural in $2$d than in $3$d.

%Despite these are all equivalent, one might select a different formulation depending on the problem at hand.

\section{Local well-posedness}\label{secLoc}

%As already pointed out above,
Because of the complicated nature of the equations, the development of a basic local wellposedness theory
(existence and uniqueness of smooth solutions for the Cauchy problem)
has %historically
proved to be highly non-trivial.
%This structure is somewhat hard to capture because of the moving domain and the quasilinear nature of the problem.
%Historically, this has made the task of establishing local wellposedness
%(existence and uniqueness of smooth solutions for the Cauchy problem) non-trivial.
Early results on the local wellposedness of the water waves system include those by Nalimov \cite{Nalimov}, Shinbrot \cite{Shin},
%Kano-Nishida \cite{KN} %-> Check this
Yosihara \cite{Yosi}, and Craig \cite{CraigLim}.
All these results deal with small perturbations of a flat interface for which the Rayleigh-Taylor
sign condition \eqref{RT} always holds.
%Other early works on the qualitative properties of the system are those by \cite{EbinMarsden,Shnirelman,Ebin1,Ebin2,BHL,DyaZak}.
%Ill-posedness when \eqref{RT} is violated has been investigated for example in \cite{BHL} and \cite{Ebin2}.
It was first observed by Wu \cite{Wu2} that, in the irrotational case, \eqref{RT} holds without smallness assumptions,
and that local-in-time solutions can be constructed with initial data of arbitrary size in Sobolev spaces \cite{Wu1,Wu2}.
Following the breakthrough of Wu,
the question of local wellposedness of the water waves and free boundary Euler equations has been addressed by several authors.
Christodoulou--Lindblad \cite{CL} and Lindblad \cite{Lindblad} considered the problem with vorticity,
Beyer--Gunther \cite{BG} took into account the effects of surface tension,
and Lannes \cite{Lannes} treated the case of non-trivial bottom topography.
The works by Ambrose-Masmoudi \cite{AM}, Coutand--Shkoller \cite{CS2}, and Shatah--Zeng \cite{ShZ1}
extended these results to more general scenarios with vorticity and surface tension,
including two-fluids systems \cite{CCS1,ShZ2,ShZ3} where surface tension is necessary for wellposedness.
Other important papers that include surface tension and/or low regularity analysis are those by
Christianson--Hur--Staffilani \cite{CHS}, and Alazard--Burq--Zuily \cite{ABZ1,ABZ2}.
%More references on the local well-posedness theory are given in \ref{RefsLoc}.
See also \cite{Schweizer,ZhangZhang,MingZhang,KinseyWu,WuAC,KuTuVi,DPN1}.
%Schweizer \cite{Schweizer}, Zhang--Zhang \cite{ZhangZhang}, Ming--Zhang \cite{MingZhang},  %Skip these?
%and the recent low regularity results by De Poyferr\'e--Nguyen \cite{DPN1}.
%and the recent works of Wu et al. \cite{KinseyWu,WuAC} on the problem of interfaces with angled crests.

Thanks to all the contributions mentioned above %and in \ref{RefsLoc} %over the past 20 years,
the local well-posedness theory is presently well-understood in a variety of different scenarios.
In short, one can say that for sufficiently nice initial configurations, it is possible to find classical smooth solutions on a small time interval.
See %section \ref{secLoc} for a discussion on the local existence theory, and
Theorem \ref{theoLoc1} for a typical result in the case of irrotational flows.

To explain some aspects of the local well-posedness theory we follow the approach in Lannes \cite{Lannes},
Alazard--M\'etivier \cite{AlMet1} and the series of works by Alazard--Burq--Zuily \cite{ABZ1}-\cite{ABZ5},
based on the use of paradifferential calculus.\footnote{See Bony \cite{Bony},
or the books of M\'etivier \cite{Metivier}, Taylor \cite{TaylorBook} for the general theory of paradifferential operators.}
%More precisely we will do this  following the work of the authors with Deng and Pausader \cite{DIPP1}, using the Weyl quantization.
We choose this path mostly because it is a good starting point for the long-time regularity theory, which we discuss in the next section.

\subsection{Paradifferential analysis}
The main objective of the paradifferential analysis of the water waves system is to formulate the
Hamiltonian system \eqref{WWE}-\eqref{defG0} for the unknowns $h$ and $\phi$
%which is a system of two scalar equations for the scalar unknowns $h$, the ``height''
%, that is the displacement of the free surface from the rest position
%and $\phi$, the velocity potential restricted to the interface.
in terms of new unknowns $h$ and $\omega$, so that the quasilinear structure of the system is apparent.
In other words, one wants to identify the terms that are responsible for the loss of derivatives
in the nonlinearity, and write the equations in a convenient form,
so that it is possible to obtain a priori energy estimates by a relatively simple procedure.
A key point is to achieve a good understanding of the Dirichlet--Neumann operator $G(h)\phi$ in \eqref{defG0}.

\subsubsection{Elements of paradifferential calculus}
Given a symbol $a=a(x,\zeta)$, $x,\zeta\in\R^d$, and a function $f \in L^2(\R^d)$,
we define the paradifferential operator $T_a f$ according to\footnote{This is the so-called Weyl quantization, which is used in \cite{DIPP},
and is particularly convenient when dealing with self-adjoint operators.
Other choices are possible to define paradifferential operators,
such as the Kohn--Nirenberg quantization used in \cite{ABZ1}-\cite{ABZ5} and \cite{ADa,ADb}.}

\begin{align}
\label{Tsigmaf}
\mathcal{F} \big( T_{a}f \big) (\xi) = \frac{1}{(2\pi)^d}\int_{\R^2}\chi\Big(\tfrac{\vert\xi-\eta\vert}{\vert\xi+\eta\vert}\Big)
  \widetilde{a}\big(\xi-\eta,\tfrac{\xi+\eta}{2}\big) \widehat{f}(\eta)\, d\eta,
\end{align}
where $\widehat{g}$ is the Fourier transform of $g$, $\widetilde{a}$ denotes the Fourier transform of $a$ in the first coordinate and
$\chi: [0,\infty) \mapsto [0,1]$ %=\varphi_{\leq -20}$.
is a smooth function supported in $[0,2^{-20}]$ and equal to $1$ on $[0,2^{-21}]$.

Notice that when $a = a(\zeta)$, \eqref{Tsigmaf} reduces to a standard Fourier multiplier. %and is a differential operator for polynomial $a$.
If instead $a = a(x)$, then $T_{a}f$ is the product of $a$ and $f$ where the frequencies of $a$ are restricted
to have size comparable or smaller than the frequencies of $f$. In particular, $T_af$ has the same regularity of $f$.
Moreover, one has the basic paradifferential decomposition of Bony \cite{Bony}:
\begin{equation*}
fg = T_f g + T_g f + \mbox{smoother terms}. %\mathcal{H}(f,g),
\end{equation*}
%where $\mathcal{H}(f,g)$ is term which has almost the regularity of $f$ and $g$ combined.
%As an example, if $a = a_1(x)a_2(\zeta)$, then $T_af$ is differential operator with non-constant coefficients that are smooth.

There are various ways to measure symbols $a=a(x,\zeta)$ and the norms of the associated paradifferential operators.
Without going into technical details one should think of the homogeneity of
$a$ with respect to the variable $\zeta$, when $|\zeta| \geq 1$,
as representing the order of the paradifferential operator, i.e., the number of derivatives acting on the function $f$,
while the dependence on $x$ is somewhat less relevant as far as regularity is concerned (which is what one cares about
to establish local existence).
A simple choice of a norm for symbols is\footnote{%In the application to evolution PDEs such as water waves, the symbols
%also depend on time $t$, but this does not play any relevant role in the local theory.
Different choices can be made depending on the specific situation at hand.
In particular, much more complicated norms have to be used when dealing with long-term regularity problems where
the dependence of the symbols on the time variable $t$ plays a crucial role;
see for example the decorated norms in Appendix A of \cite{DIPP}
where order, multiplicity, and regularity are all tracked.}
\begin{align}
{\Vert a \Vert}_{\mathcal{M}^{l}_{r,q}}
 := \sup_{\vert \alpha \vert + \vert \beta \vert \leq r} \sup_{\zeta \in \R^d} \,
%\langle t\rangle^{m(5/6-20\delta^2)+16\delta^2}
  \langle \zeta \rangle^{-l} \|\, \vert\zeta\vert^{\vert\beta\vert}
  \partial^\beta_\zeta \partial_x^\alpha a \Vert_{L^q_x(\R^d)},
\end{align}
where $l \in \R$ is the order of the symbol (and of the associated paradifferential operator), $r \in \mathbb{N}$
measures the smoothness in $x$, and $q \in [1,\infty]$ the integrability.
%Class $\Gamma^l_r$ (corresponding ot $q=\infty$) instead?

As in the case of differential operators, it is possible to establish several
algebra properties for suitable classes of paradifferential operators.
In particular, one has the mapping property
\begin{align}
a : H^s \mapsto H^{s-l}, \qquad \mbox{for} \quad \, a \in \mathcal{M}^{l}_{r,\infty},
\end{align}
and the formulas
\begin{align}
T_a T_b \approx T_{ab}, \qquad [T_a, T_b] \approx iT_{\{a,b\}}, \qquad (T_a)^\ast=T_{\bar{a}},
\end{align}
where $\approx$ denotes an equality up to terms %operators
of lower order, $[\cdot,\cdot]$ is the commutator, and $\{ a, b\} := \nabla_x a \nabla_\zeta b - \nabla_\zeta a \nabla_x b$ is the Poisson bracket.

\subsubsection{The ``good unknown'' and the Dirichlet-Neumann operator}
To describe the action of the Dirichlet-Neumann (DN) operator one introduces
\begin{align}
\label{defBV}
B := \frac{G(h)\phi+\nabla_xh\cdot\nabla_x\phi}{1+\vert\nabla h\vert^2},
  \qquad V := \nabla_x \phi - B\nabla_x h,\qquad\omega := \phi - T_B h.
\end{align}
Here $(V,B) \in \R^{n-1}\times \R$ is the restriction to the interface of the velocity field $v$ of the fluid,
and function $\omega$ is the so-called ``good-unknown'' of Alinhac \cite{AlinhacParacom,Alinhac,AlMet1}. %% Good ref here?
The origin of $\omega$ is in the paracomposition formula
$f \circ g \approx T_{f^\prime \circ g} g$, %+ \mbox{smoother terms}
which holds when $g$ is rougher than $f$. %as can be seen by differentiating the composition
%Indeed, to understand the DN operator one needs to solve the elliptic problem
%$\Delta \Phi = 0$ in $\Omega_t$ with $\Phi|_{S_t} = \phi$.
%It is then convenient to change variables $(x,y) \mapsto (x,z) = (x, y-h(t,x))$ mapping $\Omega_t$ %to the lower half-plane $H = \{z<0\}$,
%and solve the elliptic problem, now with (relatively rough) coefficients depending on $h$, in $H$.
%Since the change of variables %, or equivalently the coefficients in the new elliptic equation,
%is not regular enough, one then applies \eqref{parac} to $f = \Phi$ and $g = z+h(t,x)$:
%\begin{align}
%\Phi(x,z+h(t,x)) = T_{\partial_y\Phi (x,z+h(t,x))} \big(z+h(t,x)\big) + \mbox{smoother terms}.
%\end{align}
As a result, the variable $\omega$ in \eqref{defBV} has better smoothness properties than $\phi$,
when $h$ has limited regularity.
One of the most important outcomes of the paradifferential analysis in the context of water waves is the following key formula for the DN operator:
\begin{proposition}[Paralinearization of the DN operator]%\cite{AlMet1,ABZ1}
Let $G(h)$ be the operator defined in \eqref{defG0}, with \eqref{Phi}-\eqref{phi}, and let $\omega,B,V$ be given by \eqref{defBV}.
Then we have
\begin{align}
\label{DNformula}
G(h)\phi = T_{\lambda_{DN}}\omega - \div (T_V h) + G_2, %+ \e_1^3\mathcal{O}_{3,3/2},
\end{align}
where $G_2$ denotes smoother terms which are also quadratic in the unknowns,
and the symbol of the operator is given by
\begin{align}
\label{deflambda}
\begin{split}
\lambda_{DN} &:= \lambda^{(1)}+\lambda^{(0)},
\\
\lambda^{(1)}(x,\zeta) &:= \sqrt{(1+\vert\nabla h\vert^2)\vert\zeta\vert^2-(\zeta\cdot\nabla h)^2},
\\
\lambda^{(0)}(x,\zeta) &:= \Big( \frac{(1+|\nabla h|^2)^2}{2\lambda^{(1)}}
  \Big\{ \frac{\lambda^{(1)}}{1+|\nabla h|^2}, \frac{\zeta \cdot \nabla h}{1+|\nabla h|^2} \Big\} + \frac{1}{2}\Delta h \Big).
  %\varphi_{\geq 0}(\zeta).
\end{split}
\end{align}
\end{proposition}

%Comments on the formulas
This proposition shows the relevance of the good unknown in the fact that the main action of the DN map
can be expressed in terms of a paradifferential operator acting on $\omega$, plus a simple transport-like term $\div (T_V h)$.

%\comment{Put some classical and related references: Calderon for principal (highest order) symbol, Iooss--Plotnikov\dots}

\subsubsection{%Paralinearization,
Diagonalization and energy estimates}
Using \eqref{DNformula} and standard paralinearization arguments, one can reduce \eqref{WWE} to the following system:
\begin{align}
\label{WWpara}
\left\{
\begin{array}{l}
\partial_t h = T_{\lambda_{DN}} \omega - \div( T_V h) + G_2 %+ \varep_1^3O_{3,1},
\\ \\
\partial_t \omega = - g h  - T_\ell h - T_V \nabla\omega + \Omega_2 %+ \varep_1^3O_{3,1},
\end{array}
\right.
\end{align}
where
\begin{align}
\label{defell}
\begin{split}
\ell(x,\zeta) & := L_{ij}(x) \zeta_i \zeta_j - (g|\nabla| + \sigma|\nabla|^3)h,
  \qquad L_{ij} := \frac{\sigma }{\sqrt{1+|\nabla h|^2}} \Big(\delta_{ij} - \frac{\partial_i h \partial_j h}{1 + |\nabla h|^2} \Big),
\end{split}
\end{align}
and $G_2$ and $\Omega_2$ denote smoother terms which are quadratic in the unknowns. One can then arrive at the following result:

%At this stage it is convenient to symmetrize this system by introducing the diagonal variable

\begin{proposition}[Diagonalization and a priori energy estimates]
Let $(h,\omega)$ be solutions of \eqref{WWpara}-\eqref{defell}, and recall the notation \eqref{Tsigmaf} and \eqref{defBV}.
Define the diagonal variable\footnote{The
choice of $U$ in \eqref{defU} that symmetrizes the system \eqref{WWpara} is unique at highest order,
but can be modified by adding lower order terms.
Different choices, such as the one made in \cite{DIPP}
can be important when dealing with long-term regularity problems, where the structure of the nonlinear terms plays a major role.}
\begin{align}
\label{defU}
& U :=  T_{\sqrt{g+\ell}} h + i T_{\Sigma} T_{1/\sqrt{g+\ell}} \omega,
\end{align}
where
\begin{align}
\label{defSigma}
\Sigma & := \sqrt{\lambda_{DN}(g+\ell)} 
\end{align}
is the diagonal symbol. Then, the following hold:

\setlength{\leftmargini}{1.5em}
\begin{itemize}

\item[(i)]
$U$ satisfies the equation
\begin{align}
\label{paraeq}
\begin{split}
& \partial_t U + i T_{\Sigma} U + i T_{V\cdot\zeta}U = \mathcal{N}_U
%+ \varep_1^2 \mathcal{O}^\ast_{2,1} + \varep_1^3\mathcal{O}_{3,0},
%\\
%& Q_U :=(-\frac{1}{2}T_{\sqrt{g+\ell}\, \div V}-iT_{m^\prime (g+\ell)})h
%+ (-T_{m^\prime\Sigma}+\frac{i}{2}T_{\sqrt{\lambda_{DN}}\, \div V})\omega = 0,
%\\
%&\mathcal{N}_U := \frac{i}{3}[T_{V\cdot\zeta},T_{\Sigma}]T_{\Sigma}^{-1}(2H+i\Psi)
%  = \frac{i}{6} [T_{V\cdot\zeta},T_{\Sigma}]T_{\Sigma}^{-1} \big(3U +\overline{U} \big)+\varep_1^3\mathcal{O}_{3,0},
\end{split}
\end{align}
where $\mathcal{N}_U$ denotes nonlinear terms of lower order.

\item[(ii)] For any $k\geq 0$ consider
\begin{align}
\label{localE}
E(t) := \int_{\R^d} \big| W(t,x) \big|^2 \, dx, \qquad W := T_\Sigma^k U;
\end{align}
then\footnote{The Sobolev regularity $3k/2$ in \eqref{localEest} corresponds to the case $\sigma >0$ where
the operator $\Sigma$ is of order $3/2$. The order is instead $1/2$ in absence of surface tension.}
\begin{align}
\label{localEest}
E(t) \sim {\| U(t) \|}_{H^{3k/2}}^2 \qquad \mbox{and} \qquad \frac{d}{dt} E(t) \lesssim P\big(E(t)\big) \, E(t),
\end{align}
where $P$ is a polynomial with positive coefficients.
\end{itemize}
\end{proposition}

The procedure leading to \eqref{paraeq} is the nonlinear analogue of the basic diagonalization in \eqref{WWElin2}.
For the purpose of the local existence theory one can think of
the nonlinear terms having the form $\mathcal{N}_U = U^2$.
Then, the paralinearized equation \eqref{paraeq} has a fairly simple structure which allows one to derive the energy estimates \eqref{localEest}
in a relatively straightforward fashion, for example by applying multiple times the operator $T_\Sigma$.
%, say $k$ times, and multiplying by $(T_\Sigma)^k \bar{U}$.
\eqref{localEest} gives a priori control on a short time interval on the function $U$ in the space $H^{3k/2}$,
and hence control on $(g-\sigma\Delta)h$ and $|\nabla|^{1/2}\phi$ in the same space.
Once these a priori estimates are established, local well-posedness can then be obtained by standard procedures.
In conclusion one obtains the following result:

\begin{theorem}[Local well-posedness of \eqref{WWE}]\label{theoLoc1}
Consider the system of equations \eqref{WWE}. %, or, equivalently, the system \eqref{Euler} under the assumption \eqref{irr}.
Let an intial data $h_0 = h(t=0)$, $\phi_0 = \phi(t=0)$, be given so that with $\sqrt{g - \sigma \Delta} \,h_0$, $|\nabla|^{1/2}\phi \in H^s$,
for $s > \tfrac{d}{2} + 1$  large enough.
Then, there exists $T>0$ and a unique solution
\begin{align*}
\big(\sqrt{g - \sigma \Delta} \,h, |\nabla|^{1/2}\phi \big) \in C\big([-T,T], H^{s} \times H^{s} \big)
\end{align*}
of the system \eqref{WWE} with the assigned initial data.
\end{theorem}

We notice that, at the linear level,
\begin{equation}\label{linap}
\Sigma\approx \sqrt{g|\nabla| + \sigma |\nabla|^3} =: \Lambda_{g,\sigma}(\nabla),\qquad U\approx \sqrt{g - \sigma \Delta} \, h + i |\nabla|^{1/2} \omega.
\end{equation}
This is used to construct suitable models for global analysis, see \eqref{qwe4}--\eqref{qwe5} and \eqref{ModelWW}. 
\iffalse

Local-wellposedness also holds in the presence of vorticity:

\begin{theorem}[Local well-posedness of \eqref{Euler}]\label{theoLoc2}
Consider the system of equations \eqref{Euler} with $\sigma =0$ (resp. $\sigma >0)$.
Consider an initial domain $\Omega_0 \in H^s$ (resp. $\Omega_0 \in H^{s+1}$),
and an initial divergence-free velocity field $v_0 \in H^s(\Omega_0)$, for $s > \tfrac{d}{2} + 1$  large enough.
Assume that the Rayleigh-Taylor sign condition \eqref{RT} holds.

Then, there exists $T>0$ and a unique solution
\begin{align*}
(\Omega_t,v) \in C([0,T], H^{s} \times H^{s}), \qquad (\mbox{resp. $(\Omega_t,v) \in C([0,T], H^{s+1} \times H^{s} )$})
\end{align*}
of \eqref{Euler} with $\sigma=0$ (resp. $\sigma>0$), with the assigned initial data.
\end{theorem}

\fi

\subsection{Conclusion}\label{RefsLoc}
We have summarized the main ingredients in the local existence theory using Eulerian coordinates.
This is a natural description, which is also tied to the Hamiltonian nature of the equations, and is a good starting point for the global theory.
The other formulations described in subsection \ref{otherappr} can also be used
to develop the local theory, as in the references mentioned at the beginning of the section.
This includes, of course, well-posedness results in the presence of vorticity analogous to Theorem \ref{theoLoc1}.
%An important step is to obtain a precise understanding of the Dirichlet--Neumann operator. See also %\cite{Wilkening2} for computational aspects related to the Dirichlet--Neumann operator.
We also refer to recent work of Lannes on the interaction with floating structures \cite{LannesStruc}
%Shore problem still unpublished
and de Poyferr\'e \cite{DePoy} on emerging bottom.

\section{Global regularity and asymptotic behavior}\label{global}

The problem of global existence of solutions for water waves models is more challenging, and much fewer results have been obtained so far.
As in many other quasilinear problems, global regularity has been studied in a perturbative and dispersive setting.
Large initial data can lead to breakdown in finite time, see for example the papers on ``splash'' singularities \cite{CCFGG,CSSplash}.

In three dimensions (2D interfaces), the first global regularity results were proved by  Ger\-main-Masmoudi-Shatah \cite{GMS2} and Wu \cite{Wu3DWW} for the gravity problem ($g>0$, $\sigma=0$).
Global regularity in $3$D was also proved for the capillary problem ($g=0$, $\sigma>0$) by Germain-Masmoudi-Shatah \cite{GMSC} and for the full gravity-capillary problem ($g>0$, $\sigma>0$) by Deng-Ionescu-Pausader-Pusateri \cite{DIPP}. In the case of a finite flat bottom, global regularity was proved recently by Wang \cite{Wa2,Wa3,Wa4} in both the gravity and the capillary problems in $3$D.

In two dimensions (1D interfaces), the first long-time result for the water waves system \eqref{WWE} is due to Wu \cite{WuAG}, who showed almost-global existence for the gravity problem ($g>0$, $\sigma=0$). This was improved to global regularity by the authors in \cite{IoPu2} and, independently, by Alazard-Delort \cite{ADa,ADb}.
A different proof of Wu's $2$D almost global existence result was later given by Hunter-Ifrim-Tataru \cite{HIT},
and then complemented to a proof of global regularity in \cite{IT}. See also Wang \cite{Wa1} for a global regularity result for a class of small data of infinite energy.
For the capillary problem in $2$D, global regularity
was proved by the authors in \cite{IoPu4} and, independently,
by Ifrim-Tataru \cite{IT2} in the case of data satisfying an additional momentum condition.

We remark that all the global regularity results that have been proved so far require 3 basic assumptions: small data (small perturbations of the rest solution), trivial vorticity inside the fluid, and flat Euclidean geometry. More subtle properties are also important, such as the Hamiltonian structure of the equations, the rate of decay of the linearized waves, and the resonance structure of the bilinear wave interactions.

\subsection{Main ideas}\label{Introsecideas} The classical mechanism to establish global regularity for quasilinear equations has two main components:

\setlength{\leftmargini}{1.8em}
\begin{itemize}
  \item[(1)] Propagate control of high order energy functionals (Sobolev norms and weighted norms);
\smallskip
  \item[(2)] Prove dispersion and decay of the solution over time.
\end{itemize}

The interplay of these two aspects has been present since the seminal work of Klainerman \cite{K0,K1} on nonlinear wave equations and vector-fields,
Shatah \cite{shatahKGE} on $3$d Klein-Gordon and normal forms, Christodoulou-Klainerman \cite{CK} on the stability of Minkowski space-time,
and Delort \cite{DelortKGE} on $1$d Klein-Gordon equations.

In the last few years new methods have emerged in the study of global solutions of quasilinear evolutions, inspired by the advances in semilinear theory.
The basic idea is to combine the classical energy and vector-fields methods with refined analysis of the Duhamel formula, using the Fourier transform and carefully constructed ``designer norms''.
This is the essence of the ``method of space-time resonances'' of Germain-Masmoudi-Shatah \cite{GMS2, GMSC,GM} and Gustafson-Nakanishi-Tsai \cite{GNT1},
and of the work on plasma models and water waves in \cite{IP1,IP2,GIP,DIP,IoPu1,IoPu2,IoPu3,IoPu4,DIPP}.

In the rest of this section we illustrate the development of the these ideas in the setting of water waves by analyzing 3 systems, in increasing order of difficulty: gravity water waves in 3D, gravity water waves in 2D, and gravity-capillary water waves in 3D.

For the sake of exposition, in all three cases we take the following approach: we replace the full water waves systems with suitable simplified quasilinear models, and then outline the main ideas needed to analyze these models. The quasilinear models constructed here have two main properties: (1) they capture the essential difficulties of the global theory of the full systems, and (2) they are technically simpler than the full systems, mainly because they bypass all the difficulties of the local theory, such as the use of paradifferential calculus.

One should keep in mind that there are certain difficulties in transferring the global analysis from the model equations to the real water waves systems, mostly at the level of the energy estimates. Nevertheless, our simplified models are very useful to explain some of the key ideas of the global analysis, in problems that are more algebraically transparent.

\subsection{Gravity water waves in 3D}\label{sec3d}

We consider first the system \eqref{WWE} in 2D in the case $(g,\sigma)=(1,0)$. Global regularity in this case was proved in \cite{GMS2} and \cite{Wu3DWW}.
Here we follow essentially the exposition and the proof of Germain-Masmoudi-Shatah \cite{GMS2}; Wu's theorem in \cite{Wu3DWW} is essentially equivalent,
but involves slightly different hypothesis on the data and a very different proof (in Lagrangian coordinates, using also the Clifford algebra).

\begin{theorem}\label{ThmGlo1} Assume that $h_0,\phi_0:\mathbb{R}^2\to \mathbb{R}$ are small and smooth initial data, satisfying
\begin{equation}\label{qwe1}
\|U_0\|_{H^{N+1}}+\sum_{l\in\{1,2\}}\|\,|\nabla|^{1-p}(x_l U_0)\|_{L^2}+\sup_{t\in [0,\infty)}\langle t\rangle\|e^{-it\Lambda}U_0\|_{W^{4,\infty}}\leq\varep\leq\overline{\varep}
\end{equation}
where $N$ is sufficiently large, $\overline{\varep}$ is sufficiently small, $p=1/4$, $\langle t\rangle^{-1}=\sqrt{1+t^2}$, and
\begin{equation}\label{qwe2}
U_0:=h_0+i|\nabla|^{1/2}\phi_0,\qquad\Lambda:=\Lambda_{1,0}=\sqrt{|\nabla|}.
\end{equation}
Then there is a unique global solution $U=h+i|\nabla|^{1/2}\phi\in C([0,\infty):H^N(\mathbb{R}^2))$ of the initial-value problem \eqref{WWE} with $(g,\sigma)=(1,0)$. Moreover the solution satisfies the global bounds
\begin{equation}\label{qwe3}
\langle t\rangle^{-\delta}\|U(t)\|_{H^N}+\langle t\rangle^{-\delta}\|\,|\nabla|^{1-p}(x_l u(t))\|_{L^2}+\langle t\rangle\|U(t)\|_{W^{4,\infty}}+\|U(t)\|_{L^2}\lesssim\varep,
\end{equation}
for any $t\in[0,\infty)$ and $l\in\{1,2\}$, where $\delta=10^{-8}$ is a small constant and $u(t):=e^{it\Lambda}U(t)$ is the associated linear profile of the solution $U$.
\end{theorem}

We describe now some of the main ingredients of the proof. We highlight two main ideas: (1) the proof of high order energy estimates by symmetrization, and (2) the proof of dispersive estimates using the {\it{method of space-time resonances}}.

To simplify the exposition we replace the system \eqref{WWE} with the quasilinear evolution equation
\begin{equation}\label{qwe4}
(\partial_t+i\Lambda)U=\mathcal{N}.
\end{equation}
The quadratic nonlinearity $\mathcal{N}$ is defined by
\begin{equation}\label{qwe5}
\widehat{\mathcal{N}}(\xi):=\frac{1}{(2\pi)^2}\int_{\mathbb{R}^2}\varphi_{[-10,10]}(|\eta|/|\xi|)\sum_{l\in\{1,2\}}\widehat{V_l}(\xi-\eta)\widehat{\partial_lU}(\eta)\,d\eta,\qquad V:=\nabla(|\nabla|^{-1/2}\Re U).
\end{equation}
Here, and in the rest of the section, we use smooth cutoff functions defined as follows: we fix an even smooth function $\varphi: \R\to[0,1]$ supported in $[-2,2]$ and equal to $1$ in $[-1,1]$, and let
\begin{equation*}
\varphi_k(x) := \varphi(x/2^k) - \varphi(x/2^{k-1}) , \qquad \varphi_{I}(x):=\sum_{k\in\mathbb{I}}\varphi_k(x),
\end{equation*}
for any $k\in\mathbb{Z}$ and interval $I\subseteq\mathbb{R}$. We define also the  Littlewood-Paley projections $P_k$, $P_I$ as the operators induced by the Fourier multipliers $\varphi_k$, and $\varphi_I$ respectively. 
The equations \eqref{qwe4}--\eqref{qwe5} are a good substitute for the full system \eqref{WWE}, see the discussion in section \ref{secLoc}.

The proof relies on a bootstrap argument: we assume that $U\in C([0,T]:H^N(\mathbb{R}^2))$ is a solution of \eqref{qwe4}--\eqref{qwe5} satisfying the bootstrap hypothesis
\begin{equation}\label{qwe7}
\langle t\rangle^{-\delta}\|U(t)\|_{H^N}+\langle t\rangle^{-\delta}\sum_{l\in\{1,2\}}\|\,|\nabla|^{1-p}(x_l u(t))\|_{L^2}+\langle t\rangle\|U(t)\|_{W^{4,\infty}}+\|U(t)\|_{L^2}\leq\varep_1,
\end{equation}
for any $t\in[0,T]$, where $\varep_1:=\varep^{2/3}$, and we would like to prove the improved bounds
\begin{equation}\label{qwe8}
\langle t\rangle^{-\delta}\|U(t)\|_{H^N}+\langle t\rangle^{-\delta}\sum_{l\in\{1,2\}}\|\,|\nabla|^{1-p}(x_l u(t))\|_{L^2}+\langle t\rangle\|U(t)\|_{W^{4,\infty}}+\|U(t)\|_{L^2}\lesssim\varep.
\end{equation}
This suffices, by a simple continuity argument, since the stronger bounds \eqref{qwe8} hold at time $t=0$ due to the initial-data assumptions \eqref{qwe1}.

We remark that the bootstrap norms used in \eqref{qwe7} capture the main features of the 
nonlinear solution, namely smoothness, localization in space, and sharp pointwise decay matching the decay of linear gravity waves.

\subsubsection{Energy estimates} These are very simple in our model \eqref{qwe4}--\eqref{qwe5}: we define
\begin{equation}\label{qwe15}
W=W_N:=\langle\nabla\rangle^N U,\qquad \mathcal{E}_N(t):=\|\langle\nabla\rangle^N U(t)\|_{L^2}^2=\frac{1}{(2\pi)^2}\int_{\mathbb{R}^2}|\widehat{W}(\xi)|^2\,d\xi.
\end{equation}
Then we calculate, using the equation and symmetrization (or integration by parts)
\begin{equation}\label{qwe16}
\frac{d}{dt}\mathcal{E}_N=C\int_{\mathbb{R}^2\times\mathbb{R}^2}\widehat{W}(\eta)\overline{\widehat{W}(\xi)}\widehat{\Re U}(\xi-\eta)m(\xi,\eta)\,d\xi d\eta,
\end{equation}
where
\begin{equation}\label{qwe17}
m(\xi,\eta)=\frac{\langle\xi\rangle^N}{\langle\eta\rangle^N}\frac{(\xi-\eta)\cdot\eta}{|\xi-\eta|^{1/2}}\varphi_{[-10,10]}(|\eta|/|\xi|)-\frac{\langle\eta\rangle^N}{\langle\xi\rangle^N}\frac{(\xi-\eta)\cdot\xi}{|\xi-\eta|^{1/2}}\varphi_{[-10,10]}(|\xi|/|\eta|).
\end{equation}
The symmetrization in the symbol $m$ avoids the potential loss of derivative, and the identity above can be easily used to show that
\begin{equation}\label{qwe18}
\big|\mathcal{E}_N(t)-\mathcal{E}_N(0)\big|\lesssim \int_0^t\mathcal{E}_N(s)\cdot\|U(s)\|_{W^{4,\infty}}\,ds.
\end{equation}
This leads to the desired improved energy bound in \eqref{qwe8}. 
Notice how this step relies in a crucial way on the sharp pointwise decay of $\langle t\rangle^{-1}$ for $U(t)$.
%and slow growth of the high order energy $\mathcal{E}_N(t)$.

We remark that in the real water waves systems analyzed in \cite{GMS2} and \cite{Wu3DWW} the final result is similar (an energy inequality similar to \eqref{qwe18}), but the proof is substantially more complicated because of the quasilinear structure of the problem. In particular, the proof has to address all the difficulties of the local regularity theory of the water waves models.

\subsubsection{Dispersion and decay}
It remains to control the other terms in \eqref{qwe8}. The idea is to write the equation in terms of the linear {\it profile} $u(t) := e^{it\Lambda} U(t)$,
\begin{equation}\label{waz1.0}
\partial_t\widehat{u}(\xi,t)=\sum_{+,-}\int_{\mathbb{R}^2}e^{it[\Lambda(\xi)\mp \Lambda(\xi-\eta)\mp \Lambda(\eta)]}m_{\pm\pm}(\xi,\eta)\widehat{u_{\pm}}(\xi-\eta,t)\widehat{u_{\pm}}(\eta,t)\,d\eta,
\end{equation}
where $u_+:=u$, $u_-:=\overline{u}$, the sum is taken over choices of the signs $+,-$, and $m_{\pm\pm}$ are suitable smooth multipliers. In integral form this becomes
\begin{equation}\label{waz1}
\widehat{u}(\xi,t)=\widehat{u}(\xi,0)+\sum_{+,-}\int_0^t\int_{\mathbb{R}^2}e^{is[\Lambda(\xi)\mp \Lambda(\xi-\eta)\mp \Lambda(\eta)]}m_{\pm\pm}(\xi,\eta)\widehat{u_{\pm}}(\xi-\eta,s)\widehat{u_{\pm}}(\eta,s)\,d\eta ds.
\end{equation}
 One would like to estimate $u$ by integrating by parts either in $s$ or in $\eta$. According to \cite{GMS2}, the main contribution is expected to come from the set of quadratic space-time resonances (the stationary points of the integral)
\begin{equation}\label{waz3}
\mathcal{SR}:=\{(\xi,\eta):\,\Phi(\xi,\eta)=0,\,(\nabla_\eta\Phi)(\xi,\eta)=0,\,m(\xi,\eta)\neq 0\},
\end{equation}
where $m=m_{\pm\pm}$ and the phases $\Phi$ are defined by
\begin{equation}\label{waz3.1}
\Phi(\xi,\eta):=\Lambda(\xi)\mp \Lambda(\xi-\eta)\mp \Lambda(\eta).
\end{equation}

Since $\Lambda(\rho)=\sqrt{|\rho|}$, the first main observation is that the phases $\Phi$ only vanish when either $\xi=0$, or $\eta=0$, or $\xi-\eta=0$. In this case, however, the multipliers $m$ also vanish. In other words, there are no quadratic time resonances and one can use normal forms (integration by parts in time) to transform the quadratic terms into cubic terms.

Loss of derivative is not important at this stage of the argument, so one can integrate by parts in time and use \eqref{waz1.0}. It remains to estimate the contribution of cubic terms of the form
\begin{equation}\label{waz1.8}
\int_0^t\int_{\mathbb{R}^2\times\mathbb{R}^2}e^{is\widetilde{\Phi}(\xi,\eta,\sigma)}\frac{m(\xi,\eta)}{\Phi(\xi,\eta)}m'(\eta,\sigma)\widehat{u_{\pm}}(\xi-\eta,s)\widehat{u_{\pm}}(\eta-\sigma,s)\widehat{u_{\pm}}(\sigma,s)\,d\eta d\sigma ds,
\end{equation}
where $\widetilde{\Phi}(\xi,\eta,\sigma):=\Lambda(\xi)\mp \Lambda(\xi-\eta)\mp \Lambda(\eta-\sigma)\mp \Lambda(\sigma)$.

The resulting multiplier $\frac{m(\xi,\eta)}{\Phi(\xi,\eta)}m'(\eta,\sigma)$ is regular, so one can now analyze cubic integrals of this type using again the method of space-time resonances. An important algebraic observation, which is used in the analysis of the phases $\widetilde{\Phi}_{--+},\widetilde{\Phi}_{-+-},\widetilde{\Phi}_{+--}$ to control the weighted norms, is the identity
\begin{equation}\label{waz1.9}
\nabla_\xi\widetilde{\Phi}(\xi,\eta,\sigma)=0\quad\text{ if }\quad\nabla_{\eta,\sigma}\widetilde{\Phi}(\xi,\eta,\sigma)=0 \quad \text{ and }
\quad |\eta|,|\xi-\sigma|\ll |\xi|.
\end{equation}
This is a {\it{slow propagation of iterated resonances}} property; 
more subtle versions of this property are also important in the 3D gravity-capillary model described below,
see for example \eqref{waz9.1}.

The dispersive analysis in \cite{GMS2} is simplified by the fact that there are no quadratic space-time resonances in the problem. 
However, the basic idea of the method of space-time resonances, namely to identify these points and center the analysis around them, plays a crucial role in many other global regularity results on plasma models and water waves models. 
Further developments of these ideas, and much more sophisticated arguments, are used 
in the proof global regularity for the 3D gravity-capillary model, where one has to deal with a full set of quadratic space-time resonances. See subsection \ref{sec3dfull}.

\subsection{Gravity water waves in 2D}\label{sec2d} We consider now the system \eqref{WWE} in 1D in the case $(g,\sigma)=(1,0)$. Global regularity in this case was proved in \cite{WuAG,IoPu2,ADa,ADb,HIT,IT,Wa1}. The precise assumptions on the initial data (low frequencies, high frequencies, and the number of vector-fields involved) are not identical in these papers. We will follow mostly the setup in \cite{IoPu2}.

\begin{theorem}\label{ThmGlo2} Assume that $h_0,\phi_0:\mathbb{R}\to \mathbb{R}$ are small and smooth initial data, satisfying
\begin{equation}\label{qwe21}
\|U_0\|_{H^{N+2}}+\|x\partial_x U_0\|_{H^{N/2+1}}+\|U_0\|_Z\leq\varep\leq\overline{\varep}
\end{equation}
where $N$ is sufficiently large, $\overline{\varep}$ is a sufficiently small constant, $U_0=h_0+i|\nabla|^{1/2}\phi_0$, and
\begin{equation}\label{qwe21.1}
\|f\|_Z:=\big\|(|\xi|^\beta+|\xi|^{N/2+10})\widehat{f}(\xi)\big\|_{L^\infty_\xi},\qquad \beta:=1/4.
\end{equation}

(i) Then there is a unique global solution $U=h+i|\nabla|^{1/2}\phi\in C([0,\infty):H^N(\mathbb{R}))$ of the initial-value problem \eqref{WWE} with $(g,\sigma)=(1,0)$. The solution $U$ satisfies the global bounds
\begin{equation}\label{qwe23}
\langle t\rangle^{-\delta}\|U(t)\|_{H^N}+\langle t\rangle^{-\delta}\|S U(t)\|_{H^{N/2}}+\langle t\rangle^{1/2}\|U(t)\|_{W^{N/2+4,\infty}}\lesssim\varep,
\end{equation}
for any $t\in[0,\infty)$, where $S:=(1/2)t\partial_t+x\partial_x$ is the scaling vector-field and $\delta=10^{-8}$ is small.

(ii) The solution $U(t)$ undergoes modified (nonlinear) scattering as $t\to\infty$, i.e.
\begin{equation}\label{qwe23.5}
\lim_{t\to\infty} \mathcal{F}^{-1}\{e^{iG(\xi,t)}e^{it\Lambda(\xi)}\widehat{U}(\xi,t)\}=u^\ast_\infty\text{ in }H^{N/2}
\end{equation}
where $\Lambda(\xi)=\sqrt{|\xi|}$ and
\begin{equation}\label{qwe23.6}
G(\xi,t):=\frac{|\xi|^4}{\pi}\int_0^t|\widehat{U}(\xi,s)|^2\frac{ds}{s+1}.
\end{equation}
\end{theorem}

As before, we discuss two main ideas of the proof: (1) the {\it{quartic energy inequality}} which is needed to prove energy estimates, and (2) the construction of nonlinear profiles, to prove {\it{modified scattering}} and dispersion. As before, we use the simplified model
\begin{equation}\label{qwe50}
\begin{split}
&(\partial_t+i\Lambda)U=\mathcal{N},\\
&\widehat{\mathcal{N}}(\xi):=\frac{1}{2\pi}\int_{\mathbb{R}}\varphi_{[-10,10]}(|\eta|/|\xi|)\widehat{V}(\xi-\eta)\widehat{\partial_xU}(\eta)\,d\eta,\qquad V:=\partial_x(|\nabla|^{-1/2}\Re U),
\end{split}
\end{equation}
which is the $1D$ analogue of the model \eqref{qwe4}--\eqref{qwe5}. We use again a bootstrap argument, with the bootstrap hypothesis
\begin{equation}\label{qwe51}
\langle t\rangle^{-\delta}\|U(t)\|_{H^N}+\langle t\rangle^{-\delta}\|S U(t)\|_{H^{N/2}}+\langle t\rangle^{1/2}\|U(t)\|_{W^{N/2+4,\infty}}\leq\varep_1=\varep_0^{2/3},
\end{equation}
for a solution on some time interval $[0,T]$.

\subsubsection{Energy estimates} One can start proving energy estimates as in the 2D model, see \eqref{qwe15}--\eqref{qwe18}. These identities still hold, but the bound \eqref{qwe18} does not suffice to close the energy estimate, since the optimal $L^\infty$-type decay is $\langle t\rangle^{-1/2}$ in 1D, which is far from integrable.

The idea, which was introduced by Wu \cite{WuAG}, is to refine the energy method by proving instead a {\it{quartic energy inequality}} of the form
\begin{equation}\label{qwe51.1}
\big|\mathcal{E}'_N(t)-\mathcal{E}'_N(0)\big|\lesssim \int_0^t\mathcal{E}_N(s)\cdot\|U(s)\|^2_{W^{N/2+4,\infty}}\,ds,
\end{equation}
for a suitable functional $\mathcal{E}'_N(t)$ satisfying $\mathcal{E}'_N(t)\approx\mathcal{E}_N(t)\approx \|U(t)\|_{H^N}^2$. The point is to get two factors of $\|U(s)\|^2_{W^{N/2+4,\infty}}$ in the right-hand side, in order to have almost integrable decay.

In our model \eqref{qwe50}, a quartic energy inequality can be proved easily using the identities \eqref{qwe16}--\eqref{qwe17}.{\footnote{We remark, however, that the original proofs in \cite{WuAG} and \cite{IoPu2} used a different idea based on a nonlinear change of variables and normal forms.}} The idea is to write the bulk integrals in the right-hand side of \eqref{qwe17} in terms of the linear profiles $u=e^{it\Lambda}U$ and $w=e^{it\Lambda}W$ and integrate by parts in time. More precisely, the bulk term can be written as a linear combination of integrals of the form
\begin{equation}\label{qwe70}
\int_0^t\int_{\mathbb{R}\times\mathbb{R}}e^{-is(\Lambda(\eta)-\Lambda(\xi)\pm\Lambda(\xi-\eta))}\widehat{w}(\eta,s)\overline{\widehat{w}(\xi,s)}\widehat{u_{\pm}}(\xi-\eta,s)m(\xi,\eta)\,d\xi d\eta,
\end{equation}
where $m$ is the multiplier defined in \eqref{qwe17}. The key observation is that the phases $\Lambda(\eta)-\Lambda(\xi)\pm\Lambda(\xi-\eta)$ do not vanish, except when one of the frequencies vanishes. In this case, however, the multipliers $m$ vanish as well.

The profiles $w$ satisfy transport equations similar to \eqref{waz1.0}. Integration by parts in time and changes of variables show that the integrals in \eqref{qwe70} can be written as (1) sums of cubic boundary terms of the form
\begin{equation}\label{qwe71}
\int_{\mathbb{R}\times\mathbb{R}}e^{-is(\Lambda(\eta)-\Lambda(\xi)\pm\Lambda(\xi-\eta))}\widehat{w}(\eta,s)\overline{\widehat{w}(\xi,s)}\widehat{u_{\pm}}(\xi-\eta,s)\frac{m(\xi,\eta)}{\Lambda(\eta)-\Lambda(\xi)\pm\Lambda(\xi-\eta)}\,d\xi d\eta,
\end{equation}
where $s\in\{0,t\}$, and (2) sums of quartic space-time integrals of the form
\begin{equation}\label{qwe72}
\begin{split}
\int_0^t\int_{\mathbb{R}\times\mathbb{R}}&e^{-is(\pm\Lambda(-\xi)\pm\Lambda(\sigma)\pm\Lambda(\xi-\eta)\pm\Lambda(\eta-\sigma))}\widehat{w_{\pm}}(-\xi,s)\widehat{w_{\pm}}(\xi,s)\\
&\widehat{u_{\pm}}(\xi-\eta,s)\widehat{u_{\pm}}(\eta-\sigma,s)\widetilde{m}(\xi,\eta,\sigma)\,d\xi d\eta d\sigma.
\end{split}
\end{equation}
All the quartic space-time integrals contain two copies of $w$, two copies of $u$, and, most importantly, the multipliers $\widetilde{m}$ are regular and do not lose high-order derivatives (after suitable symmetrization). The desired inequality \eqref{qwe51.1} follows: the boundary cubic expressions in \eqref{qwe71} can be combined with the quadratic energies $\mathcal{E}_N$ to produce the energy functionals $\mathcal{E}'_N$, while the quartic space-time integrals can be estimated as claimed.

The vector-field norm can also be controlled in a similar way, by proving a similar quartic energy inequality of the form
\begin{equation}\label{qwe75}
\big|\mathcal{E}'_S(t)-\mathcal{E}'_S(0)\big|\lesssim \int_0^t(\mathcal{E}_S(s)+\mathcal{E}_N(s))\cdot\|U(s)\|^2_{W^{N/2+4,\infty}}\,ds,
\end{equation}
for a suitable functional $\mathcal{E}'_S(t)$ satisfying $\mathcal{E}'_S(t)\approx\mathcal{E}_S(t)\approx \|SU(t)\|_{H^{N/2}}^2$.

Quartic energy inequalities such as \eqref{qwe51.1} were proved and played a key role
in all the (almost) global regularity results for water waves in 2D. As explained above, the main ingredient for such an inequality to hold is the absence of bilinear time-resonances. However, the implementation is somewhat delicate in certain quasilinear problems, like water waves models, due to the potential loss of derivatives.
It can be done in some cases, for example either by using carefully constructed
nonlinear changes of variables (as in Wu \cite{WuAG}, see also \cite{IoPu2}), or the ``iterated energy method'' of Germain--Masmoudi \cite{GM}, or
the ``paradifferential normal form method'' of Alazard--Delort \cite{ADb}, or the ``modified energy method'' of Hunter--Ifrim--Tataru \cite{HIT}. These methods are largely interchangeable, as long as there are no significant quadratic resonances. See also \cite{Delo} and \cite{HITW} for earlier constructions proving quartic energy inequalities like \eqref{qwe51.1} in simpler models.

The calculation we present above, using integration by parts in time in Fourier variables, has similarities with the I-method of Colliander--Keel--Staffilani--Takaoka--Tao \cite{CKSTT1,CKSTT2}, which is used extensively in semilinear problems. One should also compare this with the more involved calculation used in energy estimates in the 3D gravity--capillary model described below.

\subsubsection{Modified scattering and decay} One can start again, as in the 2D case, from identities on the profile similar to \eqref{waz1.0}--\eqref{waz1}. The phases $\Phi(\xi,\eta)=\sqrt{|\xi|}\mp\sqrt{|\eta|}\mp\sqrt{|\xi-\eta|}$ do not vanish (except when one of the frequencies vanishes), so one can use again a normal form. As in \eqref{waz1.8} we have an identity of the type
\begin{equation}\label{qwe56}
\partial_t\widehat{u'}(\xi,t)=\sum_{+,-}\int_{\mathbb{R}\times\mathbb{R}}e^{it\widetilde{\Phi}(\xi,\eta,\sigma)}\widetilde{m}(\xi,\eta,\sigma)\widehat{u'_{\pm}}(\xi-\eta,t)\widehat{u'_{\pm}}(\eta-\sigma,t)\widehat{u'_{\pm}}(\sigma,t)\,d\eta d\sigma+\widehat{\mathcal{R}_{\geq 4}}(\xi,t),
\end{equation}
where $\widetilde{\Phi}_{\pm\pm\pm}(\xi,\eta,\sigma):=\Lambda(\xi)\mp \Lambda(\xi-\eta)\mp \Lambda(\eta-\sigma)\mp \Lambda(\sigma)$, $\widetilde{m}$ are regular multipliers, $u'$ is a suitable quadratic modification of $u$, and $\mathcal{R}_{\geq 4}$ is a quartic and higher order remainder.

The situation is different in dimension 1, compared to the dimension $d=2$ analyzed earlier, because of the slow rate of decay of solutions. In fact, it turns out that some of the terms in the right-hand side, which correspond to the cubic space-time resonances, cannot be integrated in time. These cubic space-time resonances appear only in the phases $\Phi_{++-}$, $\Phi_{+-+}$, $\Phi_{-++}$, and correspond to the frequencies $(\xi,\xi,-\xi)$, $(\xi,\xi,-\xi)$, and $(\xi,-\xi,\xi)$ respectively. To remove the non-integrable contribution one can define the nonlinear profiles $u^\ast(t)$ by
\begin{equation*}
\widehat{u^\ast}(\xi,t):=e^{iL(\xi,t)}\widehat{u'}(\xi,t),\qquad L(\xi,t):=C|\xi|^4\int_0^t|\widehat{u'}(\xi,s)|^2\frac{ds}{s+1},
\end{equation*}
 where $C$ is a suitable real constant. Using \eqref{qwe56}, one can now show that the renormalized profile $u^\ast(t)$ stays bounded and converges (quantitatively) in the $Z$ norm as $t\to\infty$,
\begin{equation*}
\|u^\ast(t_2)-u^\ast(t_1)\|_Z\lesssim 2^{-\delta m}
\end{equation*}
if $m\geq 0$ and $t_1\leq t_2\in[2^{m-1},2^{m+1}]$. This leads to global control of the solution and modified scattering, as claimed.

The idea of using nonlinear profiles and modified scattering to prove global regularity was introduced in the context of water waves by the authors in \cite{IoPu2} and Alazard-Delort in \cite{ADa,ADb}. Just like the quartic energy inequality described earlier, this idea played a key role in all the global regularity results for water waves in 2D.

\subsection{Gravity-capillary water waves in 3D}\label{sec3dfull} Finally, we consider the system \eqref{WWE} in 2D with $(g,\sigma)=(1,1)$, which was analyzed in \cite{DIPP}. Let $\Omega:=x_1\partial_2-x_2\partial_1$ denote the rotation vector-field on $\mathbb{R}^2$ and let $H^{a,b}_\Omega$ denote the space of functions defined by the norm
\begin{equation*}
\|f\|_{H^{a,b}_\Omega}:=\sum_{j\leq a}\|\Omega^jf\|_{H^b}.
\end{equation*}
The main result in \cite{DIPP} is the following global regularity theorem:

\begin{theorem}\label{MainTheo}
Assume that $\delta$ is sufficiently small, $N_0,N_1,N_3$ are sufficiently large and that the data $(h_0,\phi_0)$ satisfies
\begin{equation}\label{h0p0}
\| U_0\|_{H^{N_0}\cap H_\Omega^{N_1,N_3}} + \|U_0\|_Z= \e \leq \bar{\e},\qquad U_0:=(1-\Delta)^{1/2}h_0+i|\nabla|^{1/2}\phi_0,
\end{equation}
where $\bar{\e}$ is a sufficiently small constant and the 
$Z$ norm is explained below in \ref{secdisp}, see \eqref{waz6}. 
Then, there is a unique global solution $U=(1-\Delta)^{1/2}h+i|\nabla|^{1/2}\phi\in C\big([0,\infty) : H^{N_0}\big)$ of the system \eqref{WWE} , with $(h(0),\phi(0))=(h_0,\phi_0)$.
In addition
\begin{align}\label{mainconcl1}
(1+t)^{-\delta} {\|U(t)\|}_{H^{N_0}\cap H_\Omega^{N_1,N_3}}+\|e^{it\Lambda}U(t)\|_Z\lesssim \e_0,\qquad (1+t)^{5/6-\delta} {\|\mathcal{U}(t)\|}_{L^\infty}\lesssim \e,
\end{align}
for any $t\in[0,\infty)$, where $\Lambda:=\sqrt{|\nabla|+|\nabla|^3}$.
\end{theorem}

As before, we explain some of the main ideas, including the subtle construction of the $Z$ norm, using a simplified model. The problem is substantially more difficult in this case, and we consider the more specialized model
\begin{equation}
\label{ModelWW}
\begin{split}
& \partial_t U + i\Lambda U = \nabla V\cdot\nabla U + \tfrac{1}{2}\Delta V\cdot U,\qquad V := P_{[-10,10]}\Re U,\qquad U(0) = U_0.
\end{split}
\end{equation}
Notice that $V$ is real-valued, such that solutions of \eqref{ModelWW} satisfy the $L^2$ conservation law
\begin{equation}
\label{qaz2}
{\|U(t)\|}_{L^2} = {\|U_0\|}_{L^2}, \qquad -\infty < t < \infty.
\end{equation}
This conservation is a good substitute for the Hamiltonian structure of the original water wave systems. As before, we use a bootstrap argument, with the bootstrap hypothesis
\begin{equation}\label{exr1}
(1+t)^{-\delta} {\|U(t)\|}_{H^{N_0}\cap H_\Omega^{N_1,N_3}}+\|e^{it\Lambda}U(t)\|_Z\leq\e_1=\e_0^{2/3},\qquad t\in[0,T].
\end{equation}

\subsubsection{Energy estimates}\label{SimModEnergy} Let $W:=\langle\nabla\rangle^N U$, $\mathcal{E}_N(t):=\frac{1}{(2\pi)^2}\int_{\mathbb{R}^2}|\widehat{W}(\xi,t)|^2\,d\xi$, and calculate
\begin{equation}\label{qaz4}
\frac{d}{dt}\mathcal{E}_N=C\int_{\mathbb{R}^2\times\mathbb{R}^2}m(\xi,\eta)\widehat{W}(\eta)\widehat{\overline{W}}(-\xi)\widehat{\Re U}(\xi-\eta)\,d\xi d\eta,
\end{equation}
where
\begin{equation}\label{qaz5}
m(\xi,\eta):=\frac{(\xi-\eta)\cdot (\xi+\eta)}{2}\frac{(1+|\eta|^2)^N-(1+|\xi|^2)^N}{(1+|\eta|^2)^{N/2}(1+|\xi|^2)^{N/2}}\varphi_{[-10,10]}(\xi-\eta).
\end{equation}
This is similar to \eqref{qwe15}--\eqref{qwe17}. We notice that $m(\xi,\eta)$ satisfies
\begin{equation}\label{qaz6}
m(\xi,\eta)=\mathfrak{d}(\xi,\eta)m'(\xi,\eta),\quad\text{ where }\quad\mathfrak{d}(\xi,\eta):=\frac{[(\xi-\eta)\cdot(\xi+\eta)]^2}{1+|\xi+\eta|^2},\quad m'\approx 1.
\end{equation}
The {\it{depletion factor}} $\mathfrak{d}$ is important in establishing energy estimates, due to its correlation with the
modulation function $\Phi$ (see \eqref{qaz7.1} and \eqref{corre} below). The presence of this factor is related to the exact conservation law \eqref{qaz2}.

There is a key difference between the full gravity-capillary model and the 3D gravity model discussed earlier: the dispersion relation $\Lambda(\xi)=\sqrt{|\xi|+|\xi|^3}$ in \eqref{ModelWW} has stationary points
when $|\xi|=\gamma_0:=(2/\sqrt{3}-1)^{1/2}\approx 0.393$ (see Figure \ref{DispRelFigure} below). As a result, linear solutions can only have $|t|^{-5/6}$ pointwise decay, i.e.
\begin{equation*}
\|e^{it\Lambda}\phi\|_{L^\infty}\approx |t|^{-5/6},
\end{equation*}
even for Schwartz functions $\phi$ whose Fourier transforms do not vanish on the sphere $\{|\xi|=\gamma_0\}$. As a result, the identities \eqref{qaz4} cannot be used directly to prove energy estimates, as in the 3D gravity case. Moreover, quartic energy inequalities like \eqref{qwe51.1} also fail because there are large, codimension 1, sets of quadratic resonances, with no matching null structures (see Figure \ref{ResonantSet} below). 
New ideas, which we describe below, are needed to prove the energy bounds for this problem.

\begin{figure}[ht!]
\centering
\includegraphics[width=0.45\linewidth]{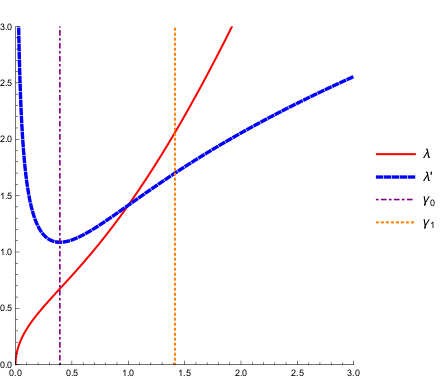} \\
\caption{\small The curves represent the dispersion relation $\lambda(r)=\sqrt{r^3+r}$ and the group velocity $\lambda'$, for $g = 1 = \sigma$. 
The frequency $\gamma_1$ corresponds to the space-time resonant sphere.
Notice that while the slower decay at $\gamma_0$ is due to some degeneracy in the linear problem, $\gamma_1$ is unremarkable from the point of view of
the linear dispersion.}
\label{DispRelFigure}
\end{figure}

{\bf{Step 1.}} We would like to estimate now the increment of $\mathcal{E}_N(t)$. We use \eqref{qaz4} and consider only the main case, 
when $|\xi|,|\eta|\approx 2^k\gg 1$, and $|\xi-\eta|$ is close to the slowly decaying frequency $\gamma_0$. 
So we need to bound space-time integrals of the form
\begin{equation*}
I:=\int_0^t\int_{\mathbb{R}^2\times\mathbb{R}^2}m(\xi,\eta)\widehat{P_kW}(\eta,s)\widehat{P_k\overline{W}}(-\xi,s)\widehat{U}(\xi-\eta,s)
\chi_{\gamma_0}(\xi-\eta)\,d\xi d\eta ds,
\end{equation*}
where $\chi_{\gamma_0}$ is a smooth cutoff function supported in the set $\{\xi:||\xi|-\gamma_0|\ll 1\}$, 
and we replaced $\Re u$ by $U$ (replacing $\Re U$ by $\overline{U}$ leads to a similar calculation). 
As before, define the linear profiles
\begin{equation}\label{qaz6.5}
u(t):=e^{it\Lambda}U(t),\qquad w(t):=e^{it\Lambda}W(t).
\end{equation}
Then decompose the integral in dyadic pieces over the size of the modulation \eqref{qaz7.1}
and over the size of the time variable. In terms of the profiles $u,w$, we need to consider the space-time integrals
\begin{equation}\label{qaz7}
\begin{split}
I_{k,m,p}:=\int_{\mathbb{R}}q_m(s)\int_{\mathbb{R}^2\times\mathbb{R}^2}e^{is\Phi(\xi,\eta)}&m(\xi,\eta)\widehat{P_kw}(\eta,s)\widehat{P_k\overline{w}}(-\xi,s)
\\
&\times\widehat{u}(\xi-\eta,s)\chi_{\gamma_0}(\xi-\eta)\varphi_p(\Phi(\xi,\eta))\,d\xi d\eta ds,
\end{split}
\end{equation}
where
\begin{align}
\label{qaz7.1}
\Phi(\xi,\eta):=\Lambda(\xi)-\Lambda(\eta)-\Lambda(\xi-\eta)
\end{align}
is the associated modulation (or phase), $q_m$ is smooth and supported in
the set $s\approx 2^m$ and $\varphi_p$ is supported in the set $\{x:|x|\approx 2^p\}$.

\begin{figure}[ht!]
\includegraphics[width=0.36\linewidth]{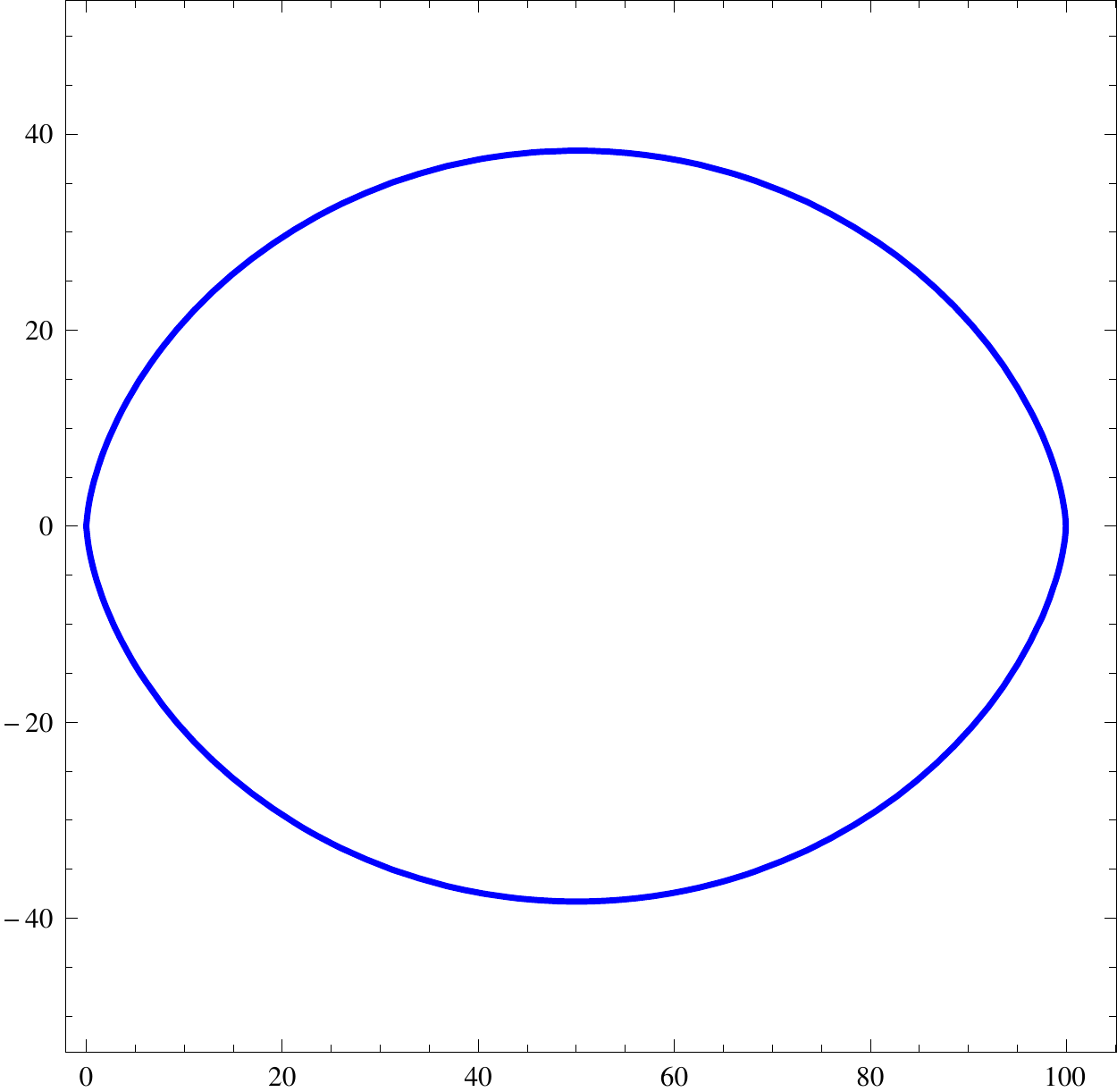}
\hskip20pt \includegraphics[width=0.403\linewidth]{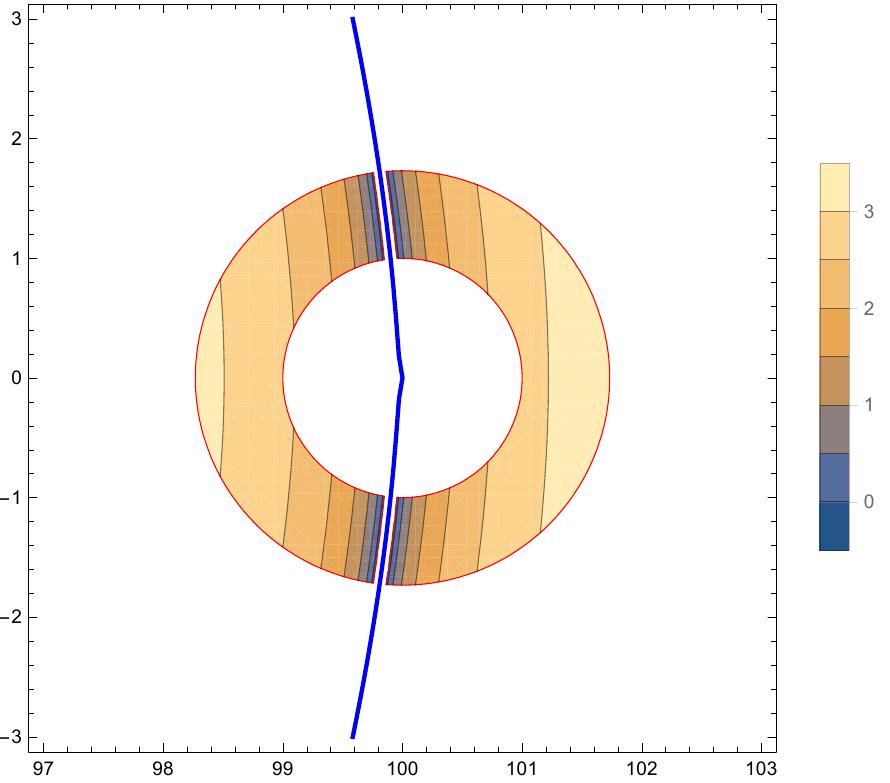} \\
\caption{\small The first picture illustrates the resonant set $\{\eta:0=\Phi(\xi,\eta)=\Lambda(\xi)-\Lambda(\eta)-\Lambda(\xi-\eta)\}$ 
for a fixed large frequency $\xi$ (in the picture $\xi=(100,0)$).
The second picture illustrates the intersection of a neighborhood of this resonant set with the set where $|\xi-\eta|$ is close to $\gamma_0$.
Note in particular that near the resonant set $\xi-\eta$ is almost perpendicular to $\xi$ (see \eqref{qaz6}, \eqref{corre}).}
\label{ResonantSet}
\end{figure}

{\bf{Step 2.}} To estimate the integrals $I_{k,m,p}$ we consider several cases depending on the relative size of $k,m,p$. 
Assume that $k,m$ are large, i.e. $2^k\gg 1, 2^m\gg 1$, which is the harder case. 
To deal with the case of small modulation, when one cannot integrate by parts in time, we need an $L^2$ bound on the Fourier integral operator
\begin{equation*}
T_{k,m,p}(f)(\xi) := \int_{\R^2} e^{is\Phi(\xi,\eta)} \varphi_k(\xi) \varphi_{\leq p}(\Phi(\xi,\eta))\chi_{\gamma_0}(\xi-\eta)
  f(\eta) \, d\eta,
\end{equation*}
where $s\approx 2^m$ is fixed. The critical bound proved in \cite{DIPP} (``the main $L^2$ lemma'') is
\begin{equation}\label{qaz9}
\|T_{k,m,p}(f)\|_{L^2}\lesssim_\eps 2^{\eps m}(2^{(3/2)(p-k/2)}+2^{p-k/2-m/3})\|f\|_{L^2},\qquad \eps>0,
\end{equation}
provided that $p-k/2\in[-0.99 m,-0.01m]$. The main gain here is the factor $3/2$ in $2^{(3/2)(p-k/2)}$ in the right-hand side 
(Schur's test would only give a factor of $1$).

The proof of \eqref{qaz9} uses a $TT^\ast$ argument, which is a standard tool to prove $L^2$ bounds for Fourier integral operators. 
This argument depends on a key nondegeneracy property of the function $\Phi$, more precisely on what we call the {\it{restricted nondegeneracy condition}}
\begin{align}\label{IntroUpsilon}
 \Upsilon(\xi,\eta) = \nabla_{\xi,\eta}^2\Phi(\xi,\eta)[\nabla_\xi^\perp\Phi(\xi,\eta), \nabla^\perp_\eta \Phi(\xi,\eta)]\neq 0\qquad\text{ if }\qquad\Phi(\xi,\eta)=0.
\end{align}
This condition, which appears to be new, %\footnote{One can compare with the case of other evolution equations, such as wave equations
%and Schr\"{o}dinger equations. In the case of wave
%equations $\Lambda(\xi)=|\xi|$, and the condition \eqref{IntroUpsilon} fails. For Schr\"{o}dinger equations $\Lambda(\xi)=|\xi|^2$, and the condition \eqref{IntroUpsilon}
%fails for the phase $\Phi(\xi,\eta)=|\xi|^2+|\eta|^2-|\xi-\eta|^2=2\xi\cdot\eta$, but holds true for other phases such as
%$\Phi(\xi,\eta)=|\xi|^2-|\eta|^2-|\xi-\eta|$.}
 can be verified explicitly in our case, when
$||\xi-\eta|-\gamma_0|\ll 1$. The function $\Upsilon$ does in fact vanish at two points on the resonant set $\{\eta:\Phi(\xi,\eta)=0\}$
(where $||\xi-\eta|-\gamma_0|\approx 2^{-k}$), but our argument can tolerate vanishing up to order $1$.

The nondegeneracy condition \eqref{IntroUpsilon} can be interpreted geometrically: the nondegeneracy of the mixed Hessian of $\Phi$ is a standard condition
that leads to optimal $L^2$ bounds on Fourier integral operators. 
In our case, however, we have the additional cutoff function $\varphi_{\leq p}(\Phi(\xi,\eta))$, so we
can only integrate by parts in the directions tangent to the level sets of $\Phi$. 
This explains the additional restriction to these directions in the
definition of $\Upsilon$ in \eqref{IntroUpsilon}.

Given the bound \eqref{qaz9}, one can control the contribution of small modulations, i.e.
\begin{equation}\label{qaz10}
 p-k/2\leq -2m/3-\epsilon m.
\end{equation}

{\bf{Step 3.}} In the high modulation case we integrate by parts in time in the formula \eqref{qaz7}. 
The main contribution is when the time derivative
hits the high frequency terms, and the resulting integral is
\begin{equation}\label{qaz12}
\begin{split}
I'_{k,m,p}:=\int_{\mathbb{R}}q_m(s)\int_{\mathbb{R}^2\times\mathbb{R}^2}e^{is\Phi(\xi,\eta)}&m(\xi,\eta)\frac{d}{ds}
\big[\widehat{P_kw}(\eta,s)\widehat{P_k\overline{w}}(-\xi,s)\big]\\
&\times\widehat{u}(\xi-\eta,s)\chi_{\gamma_0}(\xi-\eta)\frac{\varphi_p(\Phi(\xi,\eta))}{\Phi(\xi,\eta)}\,d\xi d\eta ds.
\end{split}
\end{equation}

Notice that $\partial_tw$ %satisfies the equation
%\begin{align}\label{qaz13}
%\partial_t w =\langle\nabla\rangle^Ne^{it\Lambda}\big[\nabla V\cdot\nabla U+(1/2)\Delta V\cdot U\big].
%\end{align}
%The right-hand side of \eqref{qaz13} is quadratic. 
is a quadratic expression, as in \eqref{waz1.0}, so that
%replacing $w$ by $\partial_t w$ essentially 
we gain a unit of decay (which is $|t|^{-5/6+}$), but lose a derivative. 
%This causes a problem in some range of parameters, for example when $2^{p}\approx 2^{k/2-2m/3}$,
%$1\ll 2^k\ll 2^m$, compare with \eqref{qaz10}.

%We then consider two cases: 
In the harder case when the modulation is small
we can use the depletion factor $\mathfrak{d}$ in the multiplier $m$, see \eqref{qaz6}, and the following key algebraic correlation
\begin{equation}\label{corre}
\text{ if }\qquad |\Phi(\xi,\eta)|\lesssim 1\qquad\text{ then }\qquad |m(\xi,\eta)| \lesssim 2^{-k}.
\end{equation}
See Fig. \ref{ResonantSet}. As a result, we gain one derivative in the integral $I'_{k,m,p}$, which compensates for the derivative loss.
%in \eqref{qaz13}, and the integral can be estimated again using \eqref{qaz9}.
On the other hand, when the modulation is not small, $2^p\geq 1$, then the denominator $\Phi(\xi,\eta)$ becomes a favorable factor, 
and one can %use the formula \eqref{qaz13} and
reiterate the symmetrization procedure implicit in the energy estimates. 
This avoids the loss of one derivative and gives sufficient decay to estimate $|I'_{k,m,p}|$,
and close the energy estimate.

\subsubsection{Dispersive analysis}\label{secdisp}
 The first main issue is to define an effective $Z$ norm that can be used in the bootstrap argument. 
As in \eqref{waz1}, we use the Duhamel formula, written in terms of the profile $u=u_+=e^{it\Lambda} U$, $u_-=\overline{u}$,
\begin{equation}\label{waz1.1}
\widehat{u}(\xi,t)=\widehat{u}(\xi,0)+\sum_{+,-}\int_0^t\int_{\mathbb{R}^2}e^{is[\Lambda(\xi)\mp \Lambda(\xi-\eta)\mp \Lambda(\eta)]}
m_{\pm\pm}(\xi,\eta)\widehat{u_{\pm}}(\xi-\eta,s)\widehat{u_{\pm}}(\eta,s)\,d\eta ds,
\end{equation}
where the sum is taken over choices of the signs $+,-$, and $m_{\pm\pm}$ are suitable smooth multipliers.

The idea is to estimate the function $\widehat{u}$ using the Duhamel formula \eqref{waz1}, 
by integrating by parts either in $s$ or in $\eta$. As in \eqref{waz3}--\eqref{waz3.1}, the main contribution is expected to come from the set of quadratic space-time resonances
\begin{equation}\label{waz3.2}
\mathcal{SR}:=\{(\xi,\eta):\,\Phi(\xi,\eta)=0,\,(\nabla_\eta\Phi)(\xi,\eta)=0,\,m(\xi,\eta)\neq 0\},
\end{equation}
where $m=m_{\pm\pm}$ and $\Phi(\xi,\eta)=\Lambda(\xi)\mp \Lambda(\xi-\eta)\mp \Lambda(\eta)$. 
In the gravity-capillary problem, space-time resonances are present only for the phase $\Phi(\xi,\eta)=\Lambda(\xi)-\Lambda(\xi-\eta)-\Lambda(\eta)$ and the space-time resonant set is
\begin{equation}\label{waz5}
\{(\xi,\eta)\in\mathbb{R}^2\times\mathbb{R}^2:|\xi|=\gamma_1=\sqrt{2},\,\eta=\xi/2\}.
\end{equation}
Moreover, the space-time resonant points are {\it{nondegenerate}} (according to the terminology of \cite{IP2}), 
in the sense that the Hessian of the matrix $\nabla_{\eta\eta}^2\Phi(\xi,\eta)$ is non-singular at these points.

To gain intuition, consider the first iteration of the formula \eqref{waz1}, 
i.e. assume that the functions $u_{\pm}$ in the right-hand side are Schwartz function supported at frequency $\approx 1$, 
independent of $s$. Assume that $s\approx 2^m$. 
Integration by parts in $\eta$ and $s$ shows that the main contribution comes from a small neighborhood 
of the stationary points where $|\nabla_{\eta}\Phi(\xi,\eta)|\leq 2^{-m/2+\delta m}$ and $|\Phi(\xi,\eta)|\leq 2^{-m+\delta m}$, up to negligible errors. 
Thus, the main contribution comes from space-time resonant points as in \eqref{waz3}. 
A simple calculation shows that the main contribution to the second iteration is of the type
\begin{equation*}
\widehat{u_{(2)}}(\xi)\approx c(\xi)\varphi_{\leq -m}(|\xi|-\gamma_1),
\end{equation*}
up to smaller contributions, where we have also ignored factors of $2^{\delta m}$, and $c$ is smooth.

We are now ready to describe more precisely the crucial choice of the $Z$ space.
The idea is to decompose the profile as a superposition of atoms, localized in both space and frequency,
\begin{equation*}
f={\sum}_{j,k}Q_{jk}f,\qquad Q_{jk}f=\varphi_j(x)\cdot P_kf(x).
\end{equation*}
The $Z$ norm is then defined by measuring suitably every atom.

In our case, the $Z$ space should include all Schwartz functions. It also has to include functions like $\widehat{u}(\xi)=\varphi_{\leq -m}(|\xi|-\gamma_1)$, 
due to the considerations above, for any $m$ large. 
It should measure localization in both space and frequency, and be strong enough, at least, to recover the $t^{-5/6+}$ pointwise decay. 
We define
\begin{equation}\label{waz6}
\|f\|_{Z_1}=\sup_{j,k} 2^j \cdot \Big\| \big| |\xi|-\gamma_1 \big|^{1/2}\widehat{Q_{jk}f}(\xi) \Big\|_{L^2_\xi},
\end{equation}
up to small (but important) $\delta$-corrections. Then we define the $Z$ norm by applying a suitable number of vector-fields $D$ and $\Omega$.

We emphasize that the dispersive analysis in the $Z$ norm in the gravity-capillary problem is a lot more subtle than in the earlier papers on water waves.  
To illustrate how this analysis works in our problem, we consider the contribution of the integral over $s\approx 2^m\gg 1$ in \eqref{waz1}, and assume that
the frequencies are $\approx 1$.

{\bf{Step 1.}} Start with the contribution of small modulations,
\begin{equation}\label{waz4}
\widehat{u'}(\xi):=\int_{\mathbb{R}}q_m(s)\int_{\mathbb{R}^2}\varphi_{\leq l}(\Phi(\xi,\eta))e^{is\Phi(\xi,\eta)}m_{++}(\xi,\eta)\widehat{u}(\xi-\eta,s)\widehat{u}(\eta,s)\,d\eta ds,
\end{equation}
where $l=-m+\delta m$ ($\delta$ is a small constant) and $q_m(s)$ restricts the time integral to $s\approx 2^m$, and, for simplicity, we consider only the phase $\Phi(\xi,\eta)=\Lambda(\xi)-\Lambda(\xi-\eta)-\Lambda(\eta)$. In this case the considerations above, leading to the definition of the $Z$ norm, are still relevant: one can integrate by parts in $\eta$, identify the main contributions as coming from small $2^{-m/2+\delta m}$ neighborhoods of the stationary points, and estimate these contributions in the $Z$ norm.

{\bf{Step 2.}}
Consider now the contributions of the modulations of size $2^l$, $l\geq -m+\delta m$. We start from a formula similar to \eqref{waz4} and integrate by parts in $s$. The main case is when $d/ds$ hits one of the profiles $u$. Using again the equation (see \eqref{waz1}), we have to estimate cubic expressions of the form
\begin{equation}\label{waz8}
\begin{split}
\widehat{h_{m,l}}(\xi):=\int_{\mathbb{R}}q_m(s)&\int_{\mathbb{R}^2\times\mathbb{R}^2}\frac{\varphi_{l}(\Phi(\xi,\eta))}{\Phi(\xi,\eta)}e^{is\Phi(\xi,\eta)}m_{++}(\xi,\eta)\widehat{u}(\xi-\eta,s)\\
&\times e^{is\Phi'(\eta,\sigma)}n(\eta,\sigma)\widehat{\overline{u}}(\eta-\sigma,s)\widehat{u}(\sigma,s)\,d\eta d\sigma ds,
\end{split}
\end{equation}
where $\Phi'(\eta,\sigma)=\Lambda(\eta)+\Lambda(\eta-\sigma)-\Lambda(\sigma)$. 
%Assume again that the three functions $u$ are Schwartz functions supported at frequency $\approx 1$.
We combine $\Phi$ and $\Phi'$ into the cubic phase
\begin{equation*}
\widetilde{\Phi}(\xi,\eta,\sigma):=\Phi(\xi,\eta)+\Phi'(\eta,\sigma)=\Lambda(\xi)-\Lambda(\xi-\eta)+\Lambda(\eta-\sigma)-\Lambda(\sigma).
\end{equation*}
%We need to estimate $h_{m,l}$ according to the $Z_1$ norm. Integration by parts in $\xi$ (approximate %finite speed of propagation) shows that the main contribution in $Q_{jk}h'_{m,l}$ is when $2^j\lesssim %2^m$.

%We have two main cases: if $l$ is not too small, say $l\geq -m/14$, then we use first multilinear %H\"{o}lder-type estimates, placing two of the factors $e^{is\Lambda}u$ in $L^\infty$ and one in $L^2$, %together with analysis of the stationary points of $\widetilde{\Phi}$ in $\eta$ and $\sigma$. This %suffices is most cases, except when all the variables are close to $\gamma_0$. In this case we need a %key algebraic property, of the form
%\begin{equation}\label{waz9}
%\text{ if }\quad\nabla_{\eta,\sigma}\widetilde{\Phi}(\xi,\eta,\sigma)=0\quad\text{ and %}\quad\widetilde{\Phi}(\xi,\eta,\sigma)=0\quad\text{ then %}\quad\nabla_{\xi}\widetilde{\Phi}(\xi,\eta,\sigma)=0,
%\end{equation}
%if $|\xi-\eta|,|\eta-\sigma|,|\sigma|$ are all close to $\gamma_0$.

%On the other and, if $l$ is very small, $l\leq -m/14$, then the denominator $\Phi(\xi,\eta)$ in %\eqref{waz8} is dangerous. 
The most difficult case in the dispersive analysis is when $l$ is small, say $l\leq -m/14$, and the denominator $\Phi(\xi,\eta)$ in \eqref{waz8} is dangerous.
We first restrict to suitably small neighborhoods of the stationary points of $\widetilde{\Phi}$ in $\eta$ and $\sigma$, thus to the cubic space-time resonances.
%This is the most difficult case in the dispersive analysis. 
Eventually we need to rely on one more algebraic property of the form
\begin{equation}\label{waz9.1}
\text{ if }\quad\nabla_{\eta,\sigma}\widetilde{\Phi}(\xi,\eta,\sigma)=0\quad\text{ and }\quad|\Phi(\xi,\eta)|+|\Phi'(\eta,\sigma)|\ll 1\quad\text{ then }\quad\nabla_{\xi}\widetilde{\Phi}(\xi,\eta,\sigma)=0.
\end{equation}
The point of %both \eqref{waz9} and 
\eqref{waz9.1} is that in the resonant region for the cubic integral we have $\nabla_{\xi}\widetilde{\Phi}(\xi,\eta,\sigma)=0$, so the resulting function is essentially supported when $|x|\ll 2^m$, using 
an {\it approximate finite speed of propagation} argument. This gain is reflected in the factor $2^j$ in \eqref{waz6}.

In proving control of the $Z$ norm, there are, of course, many cases to consider. 
The type of arguments presented above are typical in the proof: we decompose our profiles in space and frequency, localize to small sets in the frequency space,
keeping track in particular of the special frequencies of size $\gamma_0,\gamma_1,\gamma_1/2,2\gamma_0$, use integration by parts in $\xi$ to control the location of the output, and use multilinear H\"{o}lder-type estimates to bound $L^2$ norms. An important aspect of this analysis is that we can essentially 
assume that all profiles are almost radial and located at frequencies $\lesssim 1$,
thanks to the strong complementary control on Sobolev and weighted norms in the bootstrap hypothesis \eqref{exr1}.

{\bf{Step 3.}} The identity \eqref{waz1.1} can also be used to justify the approximate formula
\begin{equation}\label{waz7}
(\partial_t\widehat{u})(\xi,t)\approx (1/t){\sum}_j \,\, g_j(\xi)e^{it\Phi(\xi,\eta_j(\xi))}+\,\text{lower order terms},
\end{equation}
as $t\to\infty$, where $\eta_j(\xi)$ denote the stationary points where $\nabla_\eta\Phi(\xi,\eta_j(\xi))=0$. This approximate formula is consistent with the asymptotic behavior of solutions, more precisely scattering in the $Z$ norm. Qualitatively, at space-time resonances one has $\Phi(\xi,\eta_j(\xi))=0$, which leads to logarithmic growth for $\widehat{u}(\xi,t)$, while away from these space-time resonances the oscillation of $e^{it\Phi(\xi,\eta_j(\xi))}$ leads to convergence.

\subsection{Conclusions and additional references} 
To summarize, there is a small number of cases when one can construct global solutions of water waves systems, by perturbing around the trivial solutions.\footnote{In addition to the cases described earlier, there is also the capillary case $(g,\sigma)=(0,1)$,
where global solutions have been constructed in 3D in \cite{GMSC} and 2D in \cite{IT2,IoPu4}.
The proof in the capillary case follows the same path as described earlier in the gravity case,
with some additional low-frequency difficulties due to the worse dispersion relation $\Lambda(\xi)=|\xi|^{3/2}$.
See also the work of Wang \cite{Wa2}--\cite{Wa4} on finite flat bottom models in 3D.} The mechanism that leads to global solutions in all these cases is based on establishing {\it{dispersion}} and decay.

Sometimes it is possible to prove results going beyond the local theory, but not reach full global regularity. For example, starting with data of size $\e$ in a standard Sobolev space, one can sometimes get $\approx \e^{-2}$
time of existence by proving a quartic energy inequality like \eqref{qwe51.1} in cases when there are no significant quadratic resonances 
(see \cite{BIT,ITG}). 
See also the recent work of Berti--Delort \cite{BertiDelort}, where a combination of paradifferential analysis 
and ideas from KAM and normal forms theory was used to prove a significant long-time ($\approx\e^{-N}$) 
existence result for periodic 2D gravity-capillary waves (1D interface), for almost all choices of $(g,\sigma)$.
Note that these extended lifespan results do not rely on dispersion but mainly on the absence of resonances.

%not dominated by dispersive effect, and different behaviors should be expected in comparison to the full space setting described above.

In this context, a natural question to ask is if there are global or long-time regularity results for solutions with nontrivial vorticity.
We emphasize that all the global regularity results so far assume irrotationality.

\section{Formation of singularities and other topics}\label{secOther}

In this section we briefly present a few additional questions concerning the
evolution of water waves, and provide more references to other topics.

\subsection{Singularity formation}
A set of fundamental questions in pure and applied fluid dynamics concerns the study of singularities.
%While it is unreasonable at the present time to think of
While some major open problems, such as the loss of regularity and blow-up in the (rotational) Euler flow,
remain widely open, some types of ``geometric singularities'' have been studied in the context of water waves.

In \cite{CCFGL} the authors proved that a wave that is initially given as the graph of a function $h$ can overturn at a later time. 
More importantly, Castro--Cord\'oba--Gancedo--Fefferman--G\'omez-Serrano \cite{CCFGG} showed the existence of ``splash'' singularities.
A ``splash'' (resp. a ``splat'') occurs when the surface of the fluid self-intersects at a point (resp. on an arc) while retaining its smoothness.
This changes the topology of the domain, and leads to a breakdown of the chord-arc condition, that is,
the assumption (to some extent necessary for well-posedness)
\begin{align*}
\sup_{\alpha \neq \beta\in\R} \frac{\alpha - \beta}{|z(t,\alpha) - z(t,\beta)|} < \infty,
\end{align*}
where $z(t,\cdot): \R \rightarrow \R^2$ a parametrization of the interface $S_t$.
Notice that a ``splash'' singularity occurs while the parametrization of the interface and the velocity of the fluid retain their initial regularity.
The $2$d result of \cite{CCFGG} was extended to $3$ dimensions and to some other related models by Coutand--Shkoller \cite{CSSplash}.
%Similar results for the free boundary Navier-Stokes equation were proven in \cite{CCFGG-NS} and \cite{CSSplash-NS}

\begin{figure}[ht!]
{\includegraphics[width=0.6\textwidth,height=0.15\textheight]{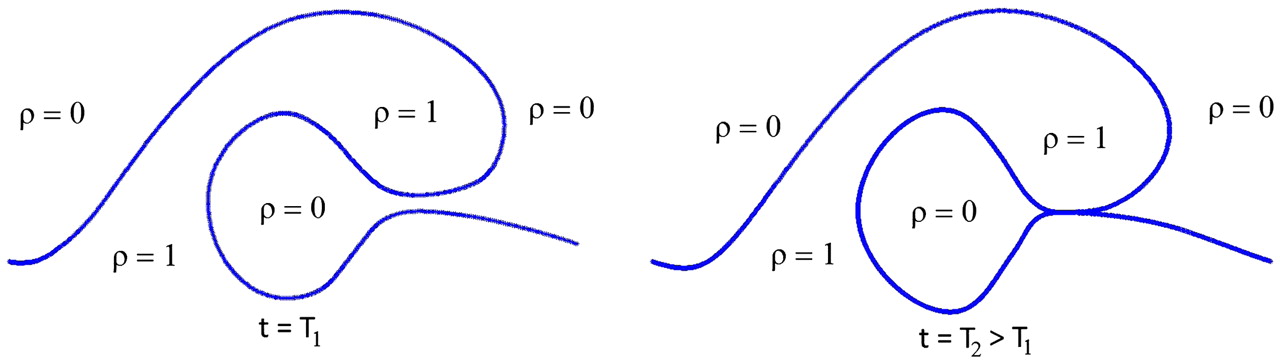}} \\
\caption{\small  Formation of a ``splash'' singularity in 2D (taken from \cite{CCFGG}).
}
  \label{fig:WWsplash}
\end{figure}

In \cite{FIL} Fefferman--Ionescu--Lie showed that a ``splash'' singularity cannot happen in the case of an interface separating two fluids:
the presence of a second fluid, with positive constant density, prevents the interface from self-intersecting.
In other words, one fluid cannot squeeze the other one if the interface and the solution are to remain smooth.
A similar result was also proven by Coutand--Shkoller in \cite{CSSnoplash}.
We refer the reader to subsection \ref{secFR2} below for some references about the well-posedness and instability for interfaces between two fluids.

%Paper by Coutand on the singularity formation in the presence of a fixed outer boundary \cite{CoutandSing}

As discussed above, a self-intersection of an interface through a fluid cannot happen for sufficiently regular solutions of the water waves equations.
However, it is plausible that a self-intersection could happen with the surface losing regularity, for example, pinching out and creating a cusp,
see Figure \ref{fig:WWsing} below.
\begin{figure}[ht!]
{\includegraphics[width=0.6\textwidth,height=0.2\textheight]{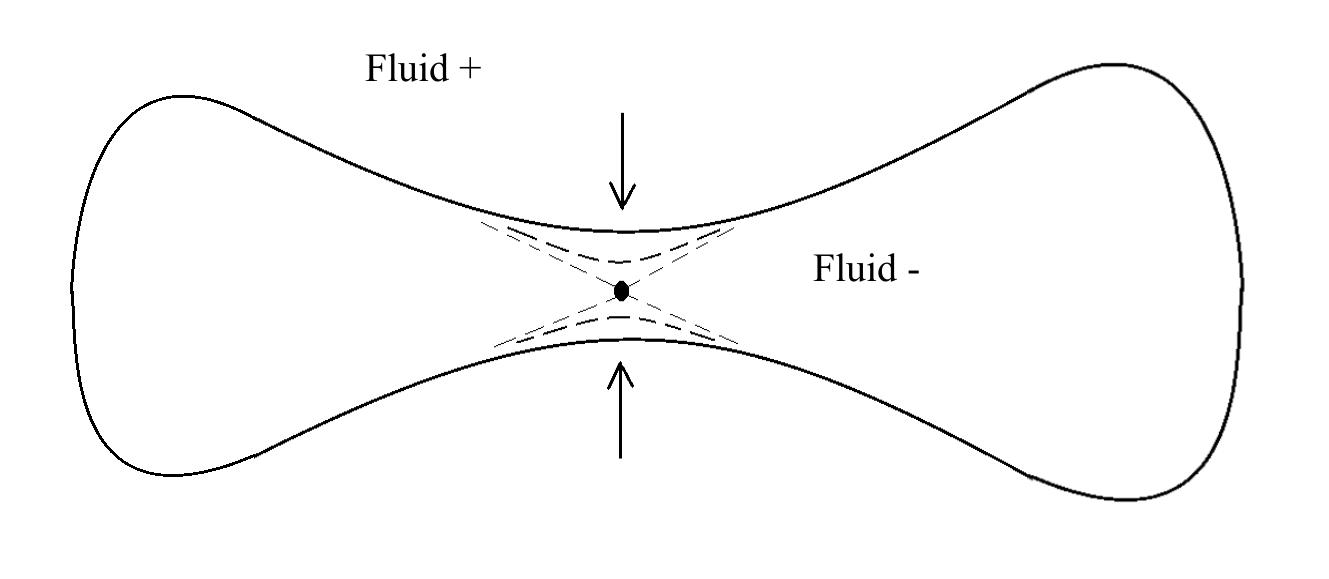}} \\ % or WWsing.jpeg
\caption{\small Possible scenario for a corner-like singularity:
(1) a locally strong velocity field pushes two points to come close together,
(2) the fluid in the middle does not have enough time to escape,
(3) a self-intersection of the interface cannot happen with a smooth boundary,
and the symmetries of the equation force the formation of a corner.
}
  \label{fig:WWsing}
\end{figure}

%%%
%%% Explain this scenario in more details?
%%%
%%% No need to...

\subsection{Fluid interfaces}\label{secFR2}
In Section \ref{secLoc} we focused our attention on the evolution problem for one fluid in vaccum.
When considering of the motion of waves on the surface of the ocean, one can think of the one-fluid model as a good first approximation
for a water-air interface; rigorous results in this direction are provided by \cite{VS,CCS2,Lannes2fluids}.
However, the motion of a free boundary between two immiscible fluids (or gases) is a more complex problem than \eqref{Euler}.

The two-fluids model is a more unstable scenario than the one-fluid model and is subject to instabilities/ill-posedness in the absence of surface tension. 
Early works on this model and the study of its instability include \cite{SS2,Ebin2,BFSS,ITT}. %,CafOre,DucRob}. % Caflish--Orellana, Duchon--Robert, Kamotski--Lebeau ?
We also mention some classical numerical works by Hou--Lowengrub--Shelley \cite{HouLowShe1,HouLowShe2}.
More recent contributions can be found in \cite{Ambrose,WuVS,ShZ3,CCS1,CCS2,VS,Lannes2fluids,AW}
to cite a few. %Put Granero-Belinchon--Shkoller?
We refer the reader to the survey of Bardos--Lannes \cite{BarLan} for more on the instability of fluids interfaces.

We remark that there are no global regularity results for any two-fluids models.

%\medskip
%\subsection{Other questions}

%Put a few open questions of interest (with some references)?
%Examples: Vorticity, Periodic problem

%\tofill

\subsection{Other questions and further references}\label{Refs}

The enormous complexity of the free boundary Euler flow and the water waves equations
has motivated %, since the very early developments of the theory, 
a large amount of research, beyond the local and global well-posedness discussed above.
This includes the construction of particular classes of solutions,
and extensive numerical activities. %from the point of view of numerical computations.
We provide here just a few additional references, including books and reviews, 
and other research papers on various topics of interest, 
and refer the reader to the cited works for more references on these topics.

%These include waves that travel horizontally at constant speed and are ``steady'' in the moving frame,
%solitary waves, %of a particular form with possibly singular interfaces, such as the famous Stokes wave %\cite{Stokes},
%standing waves that are spatially and time-periodic,  %with various periodicity properties
%waves of a particular, and possibly singular, shape such as the Stokes' wave etc\dots

%One method of investigating the properties of water waves has been
%to search for special solutions of the free-surface Euler equations for an inviscid, incompressible fluid such as
%traveling and standing waves.
%One of the earliest examples of this is due to Stokes, who in 1880 postulated that the traveling wave
%of maximum height has a crest with an internal angle of 120 degrees.

%Moreover, the complexity of the full water waves flow, and the difficulties to analyze from a rigorous %mathematical point of view
%some questions of great interest,
%have naturally motivated extensive activity from the point of view of formal asymptotics and numerical %computations.

%Given the impossibility to discuss many of the fundamental aspects in the theory of water waves 
%%besides the local and global regularity theory,
%and the great achievements of the pure and applied Mathematics community on this subject,

%\subsubsection*{Other topics}

%Particular solutions: Traveling/Steady, Standing, Solitary, Stokes...
% Simplified models
% Hamiltonian dynamics

%% Steady Water Waves:

\begin{itemize}
\item {\it Steady, solitary, extreme waves.} 
For the construction of steady water waves
 %waves that travel horizontally at constant speed and are ``steady'' in the moving frame
 we refer to Groves' survey \cite{Groves} and the paper of Constantin--Strauss \cite{ConstStrauss};
%Gerber \cite{Gerber}
the excellent review of Strauss \cite{Strausssurvey} contains an account on both the history and more modern achievements on this topic.
See also the recent paper of Constantin--Strauss--Varvaruca \cite{ConstStraussVar} for the latest developments and more references.
%Levi-Civita \cite{LeviC} and Struik \cite{Struik} on progressing waves (periodic traveling waves) in infinite depth / finite depth)

%Instability of periodic traveling waves is numerically investigated in \cite{Deconinck1}. %Deconinck--Oliveras.
%Buffoni-Groves-Sun-Wahlen \cite{BGSW} Groves?

%% Solitary waves:
The existence theory of solitary waves was developed in Friedrichs--Hyers \cite{FriedHy}, Beale \cite{Beale}, Amick--Toland \cite{AmTol1},
Amick--Kirchg\"assner \cite{AmiKir}; see also %the numerics of \cite{VBDias} and
Rousset--Tzvetkov \cite{RouTzve} for their transverse instability. %Could ask Miles for references on solitary waves etc...

%% Stokes' Wave: % Traveling wave of maximum height
The conjecture made by Stokes \cite{Stokes} %in 1848, or \cite{Stokes2} in 1880 ?
that the crest of a steady wave of maximal amplitude forms a $120^\circ$ angle has been extensively investigated.
Classical references on asymptotics are Longuet-Higgins--Fox \cite{LHFox1,LHFox2}; %on almost highest Stokes waves
a proof of the conjecture was given by Amick--Fraenkel--Toland \cite{AmFraTol};
recent numerical works dealing with the behavior of near Stokes' waves are in \cite{DyaLuKo}.
%
%% Extreme standing Waves: % regular ones below in Hamiltonian part
A counterpart of Stokes' conjecture for large standing waves was proposed by Penney--Price in the `50s
and its validity was investigated in \cite{Wilkening1}.
Recent numerical computations of $3$d standing waves can be found in Rycroft--Wilkening \cite{RycWil}.
See also \cite{LGL,ConstComm,Spielvogel,Varvaruca} for more properties of the profile of traveling gravity waves.

%% Asymptotics and Numerics:
%We also mention some classical numerical works on interfacial flows and vortex-sheets by Hou--Lowengrub--Shelley \cite{HouLowShe1,HouLowShe2},

%A recent paper by Wilkening--Vasal \cite{Wilkening2} on computational aspects of the Dirichlet-Neumann operator.

%%%% Wave Breaking ?? Not rigorous enough?

%% Hamiltonian character
\medskip
\item {\it Standing waves and Hamiltonian structure.} 
An interesting aspect of the water waves equations concerns its Hamiltonian nature, which motivates numerous questions with a strong dynamical system flavor.
An informative survey paper is Craig--Wayne \cite{CraWay}.
Works in this direction, related to small divisors and Nash-Moser techniques,
are those by Plotnikov--Toland \cite{PloTo}, Iooss--Plotnikov--Toland \cite{IoPloTo},
and Iooss--Plotnikov \cite{IoPlo} on the existence of standing waves %that is, surface water waves
that are periodic in space and time.
Quasi-periodic standing waves have been constructed by Berti--Montalto in \cite{BertiMontalto} using KAM techniques
for the first time in a quasilinear setting. See also the already mentioned work \cite{BertiDelort}
where a key role is played by the reversible, rather than Hamiltonian, structure. %Added
Aspect of the theory of normal forms and connections to integrable Hamiltonian systems are discussed in many works, see for example \cite{BenOlv, Zak2, ZakDya, CraWor, CS, CraigBirk}.

%% Simplified models
\medskip
\item {\it Approximate and asymptotics models.} Because of the physical relevance of the water waves system and the aim of better describing its complex dynamics,
many simplified models have been derived and studied in special regimes.
Two important examples include the approximation of waves in the form of wave packets by the nonlinear Schr\"odinger equation,
and the approximation of long waves in shallow water by the KdV equation.
%These include the KdV equation, the Benjamin--Ono equation, the Boussinesq and the KP equations, as well as the nonlinear Schr\"odinger equation.
We refer to \cite{CraigLim,SWgsigma,ConstLannes,ASL,TotzWuNLS}, the book \cite{LannesBook} and the surveys \cite{SWexp,Dullsurvey}
and references therein for more about reduced models and their mathematical justification.

\medskip
\item {\it Books and reviews.}
An interesting historical account on the early developments of the theory can be found in Craik \cite{Craik}.
A classical introduction to the theory of water waves is Stoker \cite{Stoker};
an introduction to water waves and related models, such as KdV and NLS, can be found in Johnson \cite{Johnson},
and a more thorough account in Sulem--Sulem \cite{SulemBook}.
The recent book of Lannes \cite{LannesBook} contains all major results on the modern well-posedness theory and on approximate models,
and the book of Constantin \cite{ConstBook} discusses many applications to oceanography.
%Both of these books posses an extensive list of references.

%Informative survey papers are those by Craig--Wayne \cite{CraWay},
%discussing various aspects from the point of view of Hamiltonian dynamics,
%Wu \cite{Wusurvey}, giving an overview on well-posedness using the complex analytical approach,
%Strauss \cite{Strausssurvey} on steady waves, and D\"ull \cite{Dullsurvey}, presenting reduced models in various regimes and results on
%and their rigorous justification.

\end{itemize}


\begin{thebibliography}{100}  %%% Remove commas after journals

\bibitem{ABZ1}%[ABZ11]
T. Alazard, N. Burq and C. Zuily.
\newblock On the water waves equations with surface tension.
\newblock {\em Duke Math. J.}, 158 (2011), no. 3, 413-499.


\bibitem{ABZ2}
T. Alazard, N. Burq and C. Zuily.
\newblock On the Cauchy problem for gravity water waves.
\newblock {\em Invent. Math.} 198 (2014), 71-163.

%% To update

\bibitem{ABZ3}%[ABZ14b]
T. Alazard, N. Burq and C. Zuily.
\newblock Strichartz estimates and the Cauchy problem for the gravity water waves equations.
\newblock Preprint. {\em arXiv:1404.4276}.


%\bibitem{ABZ4}
%T. Alazard, N. Burq and C. Zuily.
%\newblock Cauchy theory for the gravity water waves system with non localized initial data
%\newblock {\em Ann. Inst. H. Poincar\'e Anal. Non Lin\'eaire} 33 (2016), no. 2, 337-395.


\bibitem{ABZ5}
T. Alazard, N. Burq and C. Zuily.
\newblock Strichartz estimates for water waves.
\newblock {\em Ann. Sci. \'Ec. Norm. Sup\'er.}, 44 (2011), no. 5, 855-903.


\bibitem{ADa}%[AD13a]
T. Alazard and J.M. Delort.
\newblock Global solutions and asymptotic behavior for two dimensional gravity water waves.
\newblock {\em Ann. Sci. \'Ec. Norm. Sup\'er.} 48 (2015), 1149-1238.

\bibitem{ADb}%[AD13b]
T. Alazard and J.M. Delort.
\newblock Sobolev estimates for two dimensional gravity water waves
\newblock {\em  Ast\'erisque} 374 (2015) viii+241 pages.


\bibitem{AlMet1}
T. Alazard and G. M\'etivier.
\newblock Paralinearization of the Dirichlet to Neumann operator, and regularity of three-dimensional water waves.
\newblock {\em Comm. Partial Differential Equations}, 34 (2009), no. 10-12, 1632-1704.


\bibitem{AlinhacParacom}
S. Alinhac.
\newblock Paracomposition et op\'erateurs paradiff\'rentiels. %[Paracomposition and paradifferential operators]
\newblock {\em Comm. Partial Differential Equations}, 11 (1986), no. 1, 87-121.



\bibitem{Alinhac} %Correct Reference ??
S. Alinhac.
\newblock Existence d'ondes de rar\'efaction pour des syst\`emes quasi-lin\'eaires hyperboliques multidimensionnels.
\newblock {\em Comm. Partial Differential Equations}, 14 (1989), no. 2, 173-230.



\bibitem{ASL}
B. Alvarez-Samaniego and D. Lannes.
\newblock Large time existence for 3D water-waves and asymptotics.
\newblock {\em Invent. Math.} 171 (2008), no. 3, 485-541.


\bibitem{Ambrose}%[Amb03]
D.M. Ambrose.
\newblock Well-posedness of vortex sheets with surface tension.
\newblock {\em SIAM J. Math. Anal.} 35 (2003), no. 1, 211-244.

\bibitem{AM}%[AM03]
D.M. Ambrose and N. Masmoudi.
\newblock The zero surface tension limit of two-dimensional water waves.
\newblock {\em Comm. Pure Appl. Math.} 58 (2005), no. 10, 1287-1315.

\bibitem{AW}%[AM03]
D.M. Ambrose and J. Wilkening.
\newblock Computation of symmetric, time-periodic solutions of the vortex sheet with surface tension.
\newblock {\em PNAS}, vol. 107 no. 8, 3361-3366.


%%%%
%%%%

\bibitem{AmFraTol}
C. J. Amick, L. E. Fraenkel and J. F. Toland.
\newblock On the Stokes conjecture for the wave of extreme form.
\newblock {\em Acta Math.} 148 (1982), 193-214.


\bibitem{AmTol1}
C. J. Amick and J. F. Toland.
\newblock On solitary water-waves of finite amplitude.
\newblock {\em Arch. Rational Mech. Anal.} 76 (1981), no. 1, 9-95.

%\bibitem{AmTol2}
%C. J. Amick and J. F. Toland.
%\newblock The semi-analytic theory...
%\newblock {\em }

\bibitem{AmiKir}
C. J. Amick and K. Kirchg\"assner.
\newblock A theory of solitary water-waves in the presence of surface tension.
\newblock {\em Arch. Rational Mech. Anal.} 105 (1989), no. 1, 1-49.

%%%%
%%%%


\bibitem{BarLan}
C. Bardos and D. Lannes.
\newblock Mathematics for 2d interfaces.
\newblock Singularities in mechanics: formation, propagation and microscopic description, 37-67,
{\em Panor. Synth\`eses}, 38, Soc. Math. France, Paris, 2012.


\bibitem{BFSS}
C. Bardos, C. Sulem, P.-L. Sulem and U. Frisch.
\newblock Finite time analyticity for the two- and three-dimensional Kelvin-Helmholtz instability.
\newblock {\em Comm. Math. Phys.} 80 (1981), no. 4, 485-516.


\bibitem{Beale}
J. T. Beale.
\newblock The existence of solitary water waves.
\newblock {\em Comm. Pure Appl. Math.} 30 (1977), 373-389.


\bibitem{BHL}
J. T. Beale, T. Y. Hou and J. S. Lowengrub.
\newblock Growth rates for the linearized motion of fluid interfaces away from equilibrium.
\newblock {\em Comm. Pure Appl. Math.} 46 (1993), 1269-1301.


%%%

\bibitem{BenOlv}
T. B. Benjamin and P. J. Olver.
\newblock Hamiltonian structures, symmetries and conservation laws for water waves.
\newblock {\em J. Fluid Mech.} 125 (1982), 137-187.


%\bibitem{BerGer}
%F. Bernicot and P. Germain.
%\newblock Bilinear dispersive estimates via space-time resonances, part II: dimensions 2 and 3.
%\newblock {\em Arch. Ration. Mech. Anal.} 214 (2014), no. 2, 617-669.

\bibitem{BertiDelort}
M. Berti and J. M. Delort.
\newblock Almost global existence of solutions for capillarity-gravity water waves equations with periodic spatial boundary conditions.
\newblock Preprint {\em arXiv:1702.04674}.


\bibitem{BertiMontalto}
M. Berti and R. Montalto.
\newblock Quasi-periodic standing wave solutions of gravity-capillary water waves.
\newblock To appear in {\em Mem. Amer. Math. Soc}. {\em arXiv:1602.02411}.



\bibitem{BG}%[BG98]
K. Beyer and M. G\"unther.
\newblock On the Cauchy problem for a capillary drop. I. Irrotational motion.
\newblock {\em Math. Methods Appl. Sci.} 21 (1998), no. 12, 1149-1183.


\bibitem{Bony}
J. M. Bony.
\newblock Calcul symbolique et propagation des singularit\'es pour les \'equations aux d\'eriv\'ees partielles non lin\'eaires.
%[Symbolic calculus and propagation of singularities for nonlinear partial differential equations]
\newblock {\em Ann. Sci. \'Ecole Norm. Sup.} (4) 14 (1981), no. 2, 209-246.


%\bibitem{Bre}
%Y. Brenier.
%\newblock Minimal geodesics on groups of volume-preserving maps and generalized solutions of the Euler equations.
%\newblock {\em Comm. Pure Appl. Math.} 52 (1999), no .4, 411-452.

%\bibitem{CastroLannes}
%A. Castro and D. Lannes.
%\newblock  Well-posedness and shallow-water stability for a new Hamiltonian formulation of the water waves equations with vorticity.
%\newblock {em Indiana Univ. Math. J.} 64 (2015), no. 4, 1169-1270.


\bibitem{CCFGL}
A. Castro, D. C\'ordoba, C. Fefferman, F. Gancedo and M. L\'{o}pez-Fern\'{a}ndez.
\newblock Rayleigh-Taylor breakdown for the Muskat problem with applications to water waves.
\newblock {\em Ann. of Math.} 175 (2012), 909-948.


\bibitem{CCFGG}%[CCFGG13]
A. Castro, D. C\'ordoba, C. Fefferman, F. Gancedo and J. G\'{o}mez-Serrano.
\newblock Finite time singularities for the free boundary incompressible Euler equations.
\newblock {\em Ann. of Math.} 178 (2013), 1061-1134.


\bibitem{CauchyMemoirs}
A-L Cauchy.
\newblock M\'emoire sur la th'eorie de la propagation des ondes \`a la surface d'un fluide pesant d'une profondeur ind\'efinie.
\newblock {\em M\'em. Pr\'esent\'es Divers Savans Acad. R. Sci. Inst. France (Prix Acad. R. Sci., concours de 1815 et de 1816)}, I:3-312 (1827).


\bibitem{CCS1}
A. Cheng, D. Coutand and S. Shkoller.
\newblock On the motion of Vortex Sheets with surface tension.
\newblock {\em  Comm. Pure Appl. Math.} 61 (2008), no. 12, 1715-1752.

\bibitem{CCS2}%[CCS09]
A. Cheng, D. Coutand and S. Shkoller.
\newblock On the limit as the density ratio tends to zero for two perfect incompressible 3-D fluids separated by a surface of discontinuity.
\newblock {\em Comm. Partial Differential Equations}, 35 (2010), 817-845.


\bibitem{CHS}
H. Christianson, V. Hur, and G. Staffilani.
\newblock Strichartz estimates for the water-wave problem with surface tension.
\newblock {\em Comm. Partial Differential Equations} 35 (2010), no. 12, 2195-2252.


\bibitem{CL}
D. Christodoulou and H. Lindblad.
\newblock On the motion of the free surface of a liquid.
\newblock {\em Comm. Pure Appl. Math.} 53 (2000), no. 12, 1536-1602.


\bibitem{CK}
D. Christodoulou and S. Klainerman.
\newblock The global nonlinear stability of the Minkowski space.
\newblock {\em Princeton Mathematical Series} 41. Princeton University Press, Princeton, NJ, 1993. x+514 pp.


\bibitem{CKSTT1}
J. Colliander, M. Keel, G. Staffilani, H. Takaoka and T. Tao.
\newblock Sharp global well-posedness for KdV and modified KdV on $\R$ and $\mathbb{T}$.
\newblock {\em  J. Amer. Math. Soc.} 16 (2003), no. 3, 705-749.

\bibitem{CKSTT2}
J. Colliander, M. Keel, G. Staffilani, H. Takaoka and T. Tao.
\newblock Resonant decompositions and the I-method for the cubic nonlinear Schr\"odinger equation on $\R^2$.
\newblock {\em  Discrete Contin. Dyn. Syst.} 21 (2008), no. 3, 665-686.

%The idea of a cusped corner-like singularity is very promising.
%It might be worth pointing out that some numerical simulations that support this possibility were reported in the paper LGL
%but they are at best indicative since in the somewhat special context of travelling periodic
%waves they turn out to be a purely numerical artefact, cf. the considerations in

\bibitem{ConstComm}
A. Constantin.
\newblock Comment on "Steep sharp-crested gravity waves on deep water".
\newblock {\em Phys. Rev. Lett.} 93 (2004), 069402.


\bibitem{ConstBook}
A. Constantin.
\newblock Nonlinear Water Waves with applications to Wave-Current interactions and Tsunamis.
\newblock CBMS-NSF Regional Conference Series in Applied Mathematics, 81.
Society for Industrial and Applied Mathematics (SIAM), Philadelphia, PA, 2011. xii+321 pp.


\bibitem{ConstLannes}
A. Constantin and D. Lannes.
\newblock The hydrodynamical relevance of the Camassa-Holm and Degasperis-Procesi equations.
\newblock {\em Arch. Ration. Mech. Anal.} 192 (2009), no. 1, 165-186.



\bibitem{ConstStrauss}
A. Constantin and W. Strauss.
\newblock Exact steady periodic water waves with vorticity.
\newblock {\em Comm. Pure Appl. Math.} 57 (2004), no. 4, 481-527.


%In the first item of Section 4.3, perhaps it is worth adding a reference to the new developments discussed 
\bibitem{ConstStraussVar}
A. Constantin, W. Strauss and E. Varvaruca.
\newblock Global bifurcation of steady gravity water waves with critical layers, 
\newblock {\em Acta Mathematica} 217 (2016), 195-262.


\bibitem{CoutandSing}
D. Coutand.
\newblock Finite time singularity formation for moving interface Euler equations
\newblock Preprint {\em arXiv:1701.01699}.

\bibitem{CS2}%[CS07]
D. Coutand and S. Shkoller.
\newblock Well-posedness of the free-surface incompressible Euler equations with or without surface tension.
\newblock {\em  J. Amer. Math. Soc.} 20 (2007), no. 3, 829-930.


\bibitem{CSSplash}
D. Coutand and S. Shkoller.
\newblock On the finite-time splash and splat singularities for the 3-D free-surface Euler equations.
\newblock {\em  Comm. Math. Phys.} 325 (2014), 143-183.



\bibitem{CSSnoplash}
D. Coutand and S. Shkoller.
\newblock On the impossibility of finite-time splash singularities for vortex sheets.
\newblock {\em Arch. Ration. Mech. Anal.} 221 (2016), no. 2, 987-1033.




\bibitem{CraigLim}%[CR85]
W. Craig.
\newblock An existence theory for water waves and the Boussinesq and Korteweg-de Vries scaling limits.
\newblock {\em  Comm. Partial Differential Equations}, 10 (1985), no. 8, 787-1003


\bibitem{CraigBirk}%[Cra96]
W. Craig.
\newblock Birkhoff normal forms for water waves. {\em Mathematical problems in the theory of water waves (Luminy, 1995)}, 57-74.
\newblock {\em Contemp. Math.}, 200, Amer. Math. Soc., Providence, RI, 1996.



\bibitem{CraigHamWW}
W. Craig.
\newblock On the Hamiltonian for water waves.
\newblock Preprint {\em arXiv:1612.08971}.

%\bibitem{CraigND}
%W. Craig, P. Nicholls and P. David.
%\newblock Travelling two and three dimensional capillary gravity water waves.
%\newblock {\em SIAM J. Math. Anal.} 32 (2000), no. 2, 323-359.


\bibitem{CraigSS}%[Cra96]
W. Craig, U. Schanz and C. Sulem.
\newblock The modulational regime of three-dimensional water waves and the Davey-Stewartson system.
\newblock {\em Ann. Inst. H. Poincar\'e Anal. Non Lin\'eaire}, 14 (1997), no. 5, 615-667.


%\bibitem{CSgsigma}
%W. Craig and C. Sulem.
%\newblock Normal Form Transformations for Capillary-Gravity Water Waves.
%\newblock In {\em Hamiltonian Partial Differential Equations and Applications}, 73-110.
%Fields Inst. Commun., 75, {\em Fields Inst. Res. Math. Sci., Toronto, ON}, 2015.


\bibitem{CS}
W. Craig and C. Sulem.
\newblock Mapping properties of normal forms transformations for water waves.
\newblock {\em Boll. Unione Mat. Ital.} 9 (2016), no. 2, 289-318.


\bibitem{CSS}
W. Craig, C. Sulem and P.-L. Sulem.
\newblock Nonlinear modulation of gravity waves: a rigorous approach.
\newblock {\em Nonlinearity} 5 (1992), no. 2, 497-522.


%\bibitem{CraSul}
%Craig, W. and Sulem, C.
%\newblock Numerical simulation of gravity waves.
%\newblock {\em Jour. of Comp. Physics} 108 (1993), 73-83.


\bibitem{CraWay}
W. Craig and C. E. Wayne.
\newblock Mathematical aspects of surface waves on water.
% {\em Uspekhi Mat. Nauk}, 62 (2007) 95-116.
\newblock {\em Russian Math. Surveys} 62 (2007), no. 3, 453-473.


\bibitem{CraWor}
W. Craig and C. P. A. Worfolk.
\newblock An integrable normal form for water waves in infinite depth.
\newblock {\em Phys. D} 84 (1995), no. 3-4, 513-531.


\bibitem{Craik}
A.D.D. Craik.
\newblock The origins of water wave theory.
\newblock {\em Annual review of fluid mechanics}. Vol. 36, 1-28.
\newblock Annu. Rev. Fluid Mech., 36, Annual Reviews, Palo Alto, CA, 2004.


%\bibitem{Dau}
%I. Daubechies.
%\newblock Continuity statements and counterintuitive examples in connection with Weyl quantization.
%\newblock {\em J. Math. Phys.}, 24, 1453 (1983); doi: 10.1063/1.525882.


\bibitem{DePoy}
T. De Poyferr\'e.
\newblock A priori estimates for water waves with emerging bottom.
\newblock {\em arXiv:1612.04103}. Preprint.

\bibitem{DPN1}
T. de Poyferr\'e and Q.-H. Nguyen.
\newblock Strichartz estimates and local existence for the gravity-capillary waves with non-Lipschitz initial velocity.
\newblock {\em J. Differential Equations} 261 (2016), no. 1, 396-438.


%\bibitem{DPN2}
%T. De Poyferr\'e and Q.H. Nguyen.
%\newblock A paradifferential reduction for the gravity-capillary waves system at low regularity and applications.
%\newblock {\em arXiv:1508.00326}. Preprint.



%%%%%%
%%%%%% Numerics Deconinck

\bibitem{Deconinck1}
B. Deconinck and K. Oliveras.
\newblock The instability of periodic surface gravity waves.
\newblock {\em J. Fluid Mech.} 675 (2011), 141-167.

%\bibitem{Deconinck2}
%B. Deconinck and zzz
%\newblock xxx
%\newblock yyy

%%%%%%
%%%%%%

\bibitem{DelortKGE}
J.M. Delort.
\newblock Global existence and asymptotic behavior for the quasilinear Klein-Gordon equation with small data in dimension 1.
\newblock {\em Ann. Sci. \'{E}cole Norm. Sup.} 34 (2001), no. 1, 1-61.

\bibitem{Delo}
J.M. Delort.
\newblock Long-time Sobolev stability for small solutions of quasi-linear Klein-Gordon equations on the circle.
\newblock {\em Trans. Amer. Math. Soc.} 361 (2009), 4299-4365.


\bibitem{DIP}
Y. Deng, A. D. Ionescu and B. Pausader.
\newblock The Euler-Maxwell system for electrons: global solutions in 2D.
\newblock Preprint (2015).

\bibitem{DIPP}
Y. Deng, A. D. Ionescu, B. Pausader, and F. Pusateri.
\newblock Global solutions for the $3D$ gravity-capillary water waves system.
\newblock Preprint {\em arXiv:1601.05685}.

%\bibitem{DIPP1}
%Y. Deng, A. D. Ionescu, B. Pausader, and F. Pusateri.
%\newblock Global solutions for the $3D$ gravity-capillary water waves system, I: energy estimates.
%\newblock Preprint (2015).

%\bibitem{DIPP2}
%Y. Deng, A. D. Ionescu, B. Pausader. and F. Pusateri.
%\newblock Global solutions for the $3D$ gravity-capillary water waves system, II: dispersive analysis.
%\newblock Preprint (2015).


\bibitem{Dullsurvey}
W.-P. D\"ull.
\newblock On the mathematical description of water waves.
\newblock Preprint {\em arXiv:1612.06242}.



\bibitem{DKSZ} %conformal mapping variables
A. I. Dyachenko, E. A. Kuznetsov, M. Spector and V. E. Zakharov.
\newblock Analytical description of the free surface dynamics of an ideal fluid (canonical formalism and conformal mapping).
\newblock {\em Phys. Lett. A} 221 (1996), 73-79.

\bibitem{DyaLuKo} %numerics on singularity using conformal mapping variables
S. A. Dyachenko, P. M. Lushnikov and A. O. Korotkevich.
\newblock The complex singularity of a Stokes wave.
\newblock {\em JETP Lett.} 98 (2013), 767-771.



\bibitem{Ebin1}%[Ebi87]
G. Ebin.
\newblock The equations of motion of a perfect fluid with free boundary are not well posed.
\newblock {\em Comm. Partial Differential Equations} 12 (1987), no. 10, 1175-1201.


\bibitem{Ebin2}%[Ebi88]
G. Ebin.
\newblock Ill-posedness of the Raileigh-Taylor and Kelvin-Helmotz problems for for incompressible fluids.
\newblock {\em Comm. Partial Differential Equations}, 13 (1988), no. 10, 1265-1295.


\bibitem{FIL}
C. Fefferman, A. D. Ionescu and V. Lie.
\newblock On the absence of ``splash'' singularities in the case of two-fluid interfaces.
\newblock {\em Duke Math. J.} 165 (2016), no. 3, 417-462.


\bibitem{FriedHy}
K. Friedrichs and D. Hyers.
\newblock The existence of solitary waves.
\newblock {\em Comm Pure Appl. Math.} 7  (1954), 517-550.

\bibitem{Gerber}
R. Gerber.
\newblock Sur les solutions exactes des \'equations du mouvement avec surface libre d'un liquide pesant.
\newblock {\em J Math Pure Appl.} 34 (1955), 185-299.

%\bibitem{Ger}
%Germain, P.
%\newblock Global existence for coupled Klein-Gordon equations with different speeds.
%\newblock {\em Ann. Inst. Fourier} (Grenoble), 61 (2011), no. 6, 2463-2506.

\bibitem{GM}
P. Germain and N. Masmoudi.
\newblock Global existence for the Euler-Maxwell system.
\newblock {\em Ann. Sci. \'Ec. Norm. Sup\'{e}r. }(4) 47 (2014), no. 3, 469-503.


\bibitem{GMS2}
P. Germain, N. Masmoudi and J. Shatah.
\newblock Global solutions for the gravity surface water waves equation in dimension 3.
\newblock {\em Ann. of Math.}, 175 (2012), 691-754.


\bibitem{GMSC}
P. Germain, N. Masmoudi and J. Shatah.
\newblock Global solutions for capillary waves equation in dimension 3.
\newblock {\em Comm. Pure Appl. Math.} 68 (2015), no. 4, 625-687.


%\bibitem{Gerstner}
%F. Gerstner.
%\newblock Theorie der wellen.
%%\newblock {\em Abhand. Koen. Boehmischen Gesel. Wiss.}, Prague, 1802.
%\newblock {\em Ann. Phys.} 32 (1809), no. 8, 412-445.

\bibitem{Groves}
M. Groves.
\newblock Steady water waves.
\newblock {\em J. Nonl. Math. Phys.} 11 (2004), 435-460.


\bibitem{GIP}
Y. Guo, A. D. Ionescu and B. Pausader.
\newblock \newblock Global solutions of the Euler-Maxwell two-fluid system in 3D.
\newblock To appear in {\em Ann. of Math}. arXiv:1303.1060.



\bibitem{GNT1}%[GNT09]
S. Gustafson, K. Nakanishi and T. Tsai.
\newblock Scattering for the Gross-Pitaevsky equation in 3 dimensions.
\newblock {\em Commun. Contemp. Math.} 11 (2009), no. 4, 657-707.


%\bibitem{HN}%[HN98]
%Hayashi, N. and Naumkin, P.
%\newblock Asymptotics for large time of solutions to the nonlinear Schr\"{o}dinger and Hartree equations.
%\newblock {\em  Amer. J. Math.}, 120 (1998), 369-389.




\bibitem{BIT}
B. Harrop-Griffiths, M. Ifrim and D. Tataru.
\newblock Finite depth gravity water waves in holomorphic coordinates.
\newblock Preprint {\em arXiv:1607.02409}.



\bibitem{HouLowShe1}
T. Y. Hou, J. S. Lowengrub and M. Shelley.
\newblock Removing the stiffness from interfacial flows with surface tension.
\newblock {\em J. Comput. Phys.} 114 (1994), no. 2, 312-338.

\bibitem{HouLowShe2}
T. Y. Hou, J. S. Lowengrub and M. Shelley.
\newblock The long-time motion of vortex sheets with surface tension.
\newblock {\em J.  Phys. Fluids 9} (1997), no. 7, 1933-1954.

%%%%%
%%%%%



\bibitem{HITW}
J. Hunter, M. Ifrim, D. Tataru and T. Wong.
\newblock Long time solutions for a Burgers-Hilbert equation via a modified energy method.
\newblock {\em Proc. Amer. Math. Soc.} 143 (2015), 3407-3412.

\bibitem{HIT}
J. Hunter, M. Ifrim and D. Tataru.
\newblock Two dimensional water waves in holomorphic coordinates.
\newblock {\em Comm. Math. Phys.} 346 (2016), 483-552.



\bibitem{IT}
M. Ifrim and D. Tataru.
\newblock Two dimensional water waves in holomorphic coordinates II: global solutions.
\newblock {\em Bull. Soc. Math. France} 144 (2016), 369-394.


\bibitem{IT2}
M. Ifrim and D. Tataru.
\newblock The lifespan of small data solutions in two dimensional capillary water waves.
\newblock Preprint {\em arXiv:1406.5471}.

\bibitem{ITG}
M. Ifrim and D. Tataru.
\newblock Two dimensional gravity water waves with constant vorticity: I. Cubic lifespan.
\newblock Preprint {\em  arXiv:1510.07732}.

\bibitem{ITT}
T. Iguchi, N. Tanaka and A. Tani.
\newblock On the two-phase free boundary problem for two-dimensional water waves.
\newblock {\em Math. Ann.} 309 (1997), no. 2, 199-223.


\bibitem{IP1}
A. D. Ionescu and B. Pausader.
\newblock The Euler-Poisson system in 2D: global stability of the constant equilibrium solution.
\newblock {\em Int. Math. Res. Not.} (2013), no. 4, 761-826.


\bibitem{IP2}%[IPa12]
A. D. Ionescu and B. Pausader.
\newblock Global solutions of quasilinear systems of Klein-Gordon equations in 3D.
\newblock {\em J. Eur. Math. Soc.} 16 (2014), no. 11, 2355-2431.


\bibitem{IoPu1}
A. D. Ionescu and F. Pusateri.
\newblock Nonlinear fractional Schr\"{o}dinger equations in one dimension.
\newblock {\em J. Funct. Anal.} 266 (2014), 139-176.


\bibitem{IoPu2}
A. D. Ionescu and F. Pusateri.
\newblock Global solutions for the gravity water waves system in 2D.
\newblock  {\em Invent. Math.} 199 (2015), no. 3, 653-804.

\bibitem{IoPu3}
A. D. Ionescu and F. Pusateri.
\newblock Global analysis of a model for capillary water waves in 2D.
\newblock {\em Comm. Pure Appl. Math.} 69 (2016), no. 11, 2015-2071.

\bibitem{IoPu4}
A. D. Ionescu and F. Pusateri.
\newblock Global regularity for 2d water waves with surface tension.
\newblock To appear in {\em Mem. Amer. Math. Soc}. {\em arXiv:1408.4428}.




%%%%%%
%%%%%%

\bibitem{IoPlo}
G. Iooss and P. I. Plotnikov.
\newblock Small divisor problem in the theory of three-dimensional water gravity waves.
\newblock {\em Mem. Amer. Math. Soc.} 200 (2009), no. 940, viii+128 pp.
%doubly periodic travelling waves on the surface of an infinitely deep three-dimensional perfect fluid
%travelling symmetric ``diamond''

\bibitem{IoPloTo}
G. Iooss, P. I. Plotnikov and J. F. Toland.
\newblock Standing waves on an infinitely deep perfect fluid under gravity.
\newblock {\em Arch. Ration. Mech. Anal.} 177 (2005), no. 3, 367-478.

%%%%%%
%%%%%%




\bibitem{John}
F. John.
\newblock Blow-up for quasilinear wave equations in three space dimensions.
\newblock {\em Comm. on Pure Appl. Math.} 34 (1981), no. 1, 29-51.

%\bibitem{JK}%[JK84]
%John, F. and Klainerman, S.
%\newblock Almost global existence to nonlinear wave equations in three space dimensions.
%\newblock {\em Comm. Pure Appl. Math.} 37 (1984), no. 4, 443-455.

\bibitem{Johnson}%[Joh81]
R. S. Johnson.
\newblock A modern introduction to the mathematical theory of water waves.
\newblock Cambridge Texts in Applied Mathematics. {\em Cambridge University Press, Cambridge}, 1997. xiv+445 pp.



\bibitem{KN}
T. Kano and T. Nishida.
\newblock Sur les ondes de surface de l'eau avec une justification math\'ematique des \'equations des ondes en eau peu profonde.
\newblock {\em J. Math. Kyoto Univ.} 19 (1979), no. 2, 335-370.


\bibitem{KinseyWu}
R. H. Kinsey and S. Wu.
\newblock A Priori Estimates for Two-Dimensional Water Waves with Angled Crests
\newblock Preprint {\em arXiv:1406.7573}.


\bibitem{K0}%[Kla85a]
S. Klainerman.
\newblock Uniform decay estimates and the Lorentz invariance of the classical wave equation.
\newblock {\em  Comm. Pure Appl. Math.} 38 (1985), no. 3, 321-332.

%\bibitem{KKG}%[Kla85b]
%Klainerman, S.
%\newblock Global existence of small amplitude solutions to nonlinear Klein-Gordon equations in four space-time dimensions.
%\newblock {\em Comm. Pure Appl. Math.}, 38 (1985), no. 5, 631-641.


\bibitem{K1}%[Kla86]
S. Klainerman.
\newblock The null condition and global existence for systems of wave equations.
\newblock {\em Nonlinear systems of partial differential equations in applied mathematics, Part 1 (Santa Fe, N.M., 1984)}, 293-326.
\newblock Lectures in Appl. Math., 23, Amer. Math. Soc., Providence, RI, 1986.

%%%%%%

\bibitem{KuTuVi}
I. Kukavica, A. Tuffaha and V. Vicol.
\newblock On the Local Existence and Uniqueness for the 3D Euler Equation with a Free Interface.
\newblock {\em Appl. Math. Optim.} (2016), 1-29.

%%%%%%

%\bibitem{KSZ1}
%E. A. Kuznetsov, M. D. Spector and V. E. Zakharov.
%\newblock Surface singularities of ideal fluid.
%\newblock {\em Phys. Lett. A} 182 (1993), no. 4-6, 387-393.

%\bibitem{KSZ2}
%E. A. Kuznetsov, M. D. Spector and V. E. Zakharov.
%\newblock Formation of singularities on the free surface of an ideal fluid.
%\newblock {\em Phys. Rev. E} 49 (1994), no. 2, 1283-1290.

%%%%%%


%\bibitem{Lagrange}
%J-L. Lagrange.
%\newblock M\'emoire sur la th'eorie du mouvement des fluides.
%\newblock ... (1781)


\bibitem{Lannes}%[Lan05]
D. Lannes.
\newblock Well-posedness of the water waves equations.
\newblock {\em J. Amer. Math. Soc.} 18 (2005), no. 3, 605-654.

\bibitem{Lannes2fluids}
D. Lannes.
\newblock A stability criterion for two-fluid interfaces and applications.
\newblock {\em Arch. Ration. Mech. Anal.} 208 (2013), no. 2, 481-567.

\bibitem{LannesBook}
D. Lannes.
\newblock The water waves problem. Mathematical analysis and asymptotics.
\newblock Mathematical Surveys and Monographs, Vol. 188. American Mathematical Society, Providence, RI, 2013. xx+321 pp.

\bibitem{LannesStruc}
D. Lannes.
\newblock On the dynamics of floating structures.
\newblock Preprint {\em arXiv:1609.06136}.


%\bibitem{Laplace}
%Laplace, P-S. Marquis de
%\newblock Suite de r\'echerches sur plusieurs points du syst\`eme du monde
%\newblock ... (1776)

\bibitem{LeviC}
T. Levi-Civita.
\newblock D\'etermination rigoureuse des ondes permanentes d'ampleur finie.
\newblock {\em Math. Ann.} 93 (1925), no. 1, 264-314.


\bibitem{Lindblad}%[Lin05]
H. Lindblad.
\newblock Well-posedness for the motion of an incompressible liquid with free surface boundary.
\newblock {\em Ann. of Math.} 162 (2005), no. 1, 109-194.


%\bibitem{LR2}%[LR05]
%H. Lindblad and I. Rodnianski.
%\newblock Global existence for the Einstein vacuum equations in wave coordinates.
%\newblock {\em Comm. Math. Phys.} 256 (2005), no. 1, 43-110.


\bibitem{LHFox1}
M. S. Longuet-Higgins and M. J. H. Fox.
\newblock Theory of the almost-highest wave: The inner solution. %Part 1: The inner solution.
\newblock {J . Fluid Mech.} (1977), vol. 80, 721-741.
%\newblock Part 2: Matching and analytic extension. {\em J . Fluid Mech.} (1978), vol. 85, 769-786.

\bibitem{LHFox2}
M. S. Longuet-Higgins and M. J. H. Fox.
\newblock Theory of the almost-highest wave. Part 2. Matching and analytic extension
\newblock {\em J . Fluid Mech.} 85 (1978), 769-786.


\bibitem{LGL}
V. Lukomsky, I. Gandzha, and D. Lukomsky.
\newblock Steep sharp-crested gravity waves on deep water. 
\newblock {\em Phys. Rev. Lett.} 89 (2002), 164502.


\bibitem{Metivier}
G. M\'etivier.
\newblock Para-differential calculus and applications to the Cauchy problem for nonlinear systems.
\newblock Centro di Ricerca Matematica Ennio De Giorgi (CRM) Series, 5. {\em Edizioni della Normale, Pisa}, 2008. xii+140 pp.

%\bibitem{MinRouTzve}
%M. Ming, F. Rousset and N. Tzvetkov.
%\newblock Multi-solitons and Related Solutions for the Water-waves System
%\newblock {\em SIAM J. Math. Anal.}, 47 (2015), no. 1, 897-954.

\bibitem{MingZhang}
M. Ming and Z. Zhang.
\newblock Well-posedness of the water-wave problem with surface tension.
\newblock {\em J. Math. Pures Appl.} (9) 92 (2009), no. 5, 429-455.


\bibitem{Nalimov}%[Na74]
V. I. Nalimov.
\newblock The Cauchy-Poisson problem.
\newblock {\em Dinamika Splosn. Sredy Vyp.} 18 Dinamika Zidkost. so Svobod. Granicami (1974), 10-210, 254.

% %\bibitem{Olv}
%P. J. Olver.
%Conservation laws of free boundary problems and the classification of conservation laws for water waves.
%\newblock {\em Trans. Amer. Math. Soc.} 277 (1983), 353-380.


% CHECK and ADD this also in relation to Holomorphic Coordinates

%\bibitem{Ovsy}
%L. V. Ovsyannikov.
%\newblock Dynamics of a fluid.
%\newblock {\em M.A. Lavrent'ev Institute of Hydrodynamics} Sib. Branch USSR Ac. Sci. 15 (1973), 104-125.


%\bibitem{PoissonMemoirs}
%S. D. Poisson.
%\newblock M\'emoire sur la th\'eorie des ondes.
%\newblock {\em Mem. Acad. R. Sci. Inst. France}. 1816, 2nd Ser. 1:70-186.

%%%

\bibitem{PloTo}
P. I. Plotnikov and J. F. Toland.
\newblock Nash-Moser theory for standing waves. %existence of standing waves in finite depth
\newblock {\em Arch. Ration. Mech. Anal.} 159 (2001), no. 1, 1-83.


\bibitem{VS}%[Pu11]
F. Pusateri.
\newblock On the limit as the surface tension and density ratio tend to zero for the two phase Euler equation.
\newblock {\em J. Hyperbolic Differ. Equ.} 8 (2011), no. 2, 347-373.


\bibitem{RycWil}
C. H. Rycroft and J. Wilkening.
\newblock Computation of three-dimensional standing water waves.
\newblock {\em J. Comput. Phys.} 255 (2013), 612-638.


\bibitem{RouTzve} %% AmiKir for existence
F. Rousset and N. Tzvetkov.
\newblock Transverse instability of the line solitary water-waves.
\newblock {\em Invent. Math.} 184 (2011), no. 2, 257-388.


%\bibitem{SWlong} %% LONG WAVE LIMIT
%G. Schneider and C. E. Wayne.
%\newblock The long wave limit for the water wave problem. I. the case of zero surface tension.
%\newblock {\em Comm. Pure Appl. Math.} 53 (2000), no. 12, 1475-1535.


\bibitem{SWexp} %% Expository paper
G. Schneider and C. E. Wayne.
\newblock On the validity of 2D-surface water wave models.
\newblock {\em GAMM Mitt. Ges. Angew. Math. Mech.} 25 (2002), no. 1-2, 127-151.
%\newblock \hyperlink{http://math.bu.edu/people/cew/preprints/models.ps}{(Web Link)}.


\bibitem{SWgsigma}
G. Schneider and C. E. Wayne.
\newblock The rigorous approximation of long-wavelength capillary-gravity waves.
\newblock {\em Arch. Ration. Mech. Anal.} 162 (2002), no. 3, 247-285.



\bibitem{Schweizer}
B. Schweizer.
\newblock On the three-dimensional Euler equations with a free boundary subject to surface tension.
\newblock {\em Ann. Inst. H. Poincar\'e Anal. Non Lin\'eaire} 22 (2005), 753-781.


\bibitem{shatahKGE}%[Sha85]
J. Shatah.
\newblock Normal forms and quadratic nonlinear Klein-Gordon equations.
\newblock {\em   Comm. Pure Appl. Math.} 38 (1985), no. 5, 685-696.


\bibitem{ShZ1}%[ShZ08a]
J. Shatah and C. Zeng.
\newblock Geometry and a priori estimates for free boundary problems of the Euler equation.
\newblock {\em Comm. Pure Appl. Math.} 61 (2008), no. 5, 698-744.

\bibitem{ShZ2} %[SZ08b]
J. Shatah and C. Zeng.
\newblock A priori estimates for fluid interface problems.
\newblock {\em  Comm. Pure Appl. Math.} 61 (2008), no. 6, 848-876.

\bibitem{ShZ3} %[ShZ12]
J. Shatah and C. Zeng.
\newblock Local well-posedness for the fluid interface problem.
\newblock {\em Arch. Ration. Mech. Anal.} 199 (2011), no. 2, 653-705.



\bibitem{Shin}
M.  Shinbrot.
\newblock The initial value problem for surface waves under gravity. I. The simplest case.
\newblock {\em Indiana Univ. Math. J.} 25 (1976), no. 3, 281-300.


\bibitem{Sideris}
T. C. Sideris.
\newblock Formation of singularities in three-dimensional compressible fluids.
\newblock {\em Comm. Math. Phys.} 101 (1985), no. 4, 475-485.

%\bibitem{Spirn2}%[SW12]
%Spirn, D. and Wright, D.
%\newblock Linear dispersive decay estimates for the 3+1 dimensional water wave equation with surface tension.
%\newblock {\em Canad. Math. Bull.} 55 (2012), no. 1, 176-187.



%Also, it is worth pointing out that the profile of a travelling surface
%gravity water wave in irrotational flow must always be a graph (overhanging is not possible), cf. the considerations pioneered in the paper

\bibitem{Spielvogel}
E. R. Spielvogel.
\newblock A variational principle for waves of infinite depth.
\newblock {\em Arch. Rational Mech. Anal.} 39 (1970), 189-205.

%and completed in Varvaruca


\bibitem{Stoker}
J. J. Stoker.
\newblock Water waves. The mathematical theory with applications.
\newblock Wiley Classics Library. A Wiley-Interscience Publication. {\em John Wiley \& Sons, Inc., New York}, 1992. xxvi+567 pp.


%%% REF % Check right reference for the extreme wave

\bibitem{Stokes}
G. G. Stokes.
\newblock On the theory of oscillatory waves.
\newblock {\em Trans. Cambridge Philos. Soc.}, 8 (1847), 441-455.

\bibitem{Stokes2}
G. G. Stokes.
\newblock Considerations relative to the greatest height of oscillatory irrotational waves which can be propagated
without change of form (1880).
\newblock {\em Mathematical and Physical Papers}, Volume 1, pages 225-228.

%%%


\bibitem{Strausssurvey}
W. Strauss.
\newblock Steady water waves.
\newblock {Bull. Amer. Math. Soc. (N.S.)} 47 (2010), no. 4, 671-694.


\bibitem{Struik}
D. J. Struik.
\newblock D\'etermination rigoureuse des ondes irrotationelles p\'eriodiques dans un canal \`a profondeur finie.
\newblock {\em Math. Ann.} 95 (1926), no. 1, 595-634.


\bibitem{SS2}
C. Sulem and P.L. Sulem.
\newblock The well-posedness of two-dimensional ideal flow.
\newblock {\em J. M\'ec. Th\'eor. Appl.} 1983, Special Issue, 217-242.


\bibitem{SulemBook}%[SS99]
C. Sulem and P.L. Sulem.
\newblock The nonlinear Schr\"{o}dinger equation. Self-focussing and wave collapse.
\newblock {\em Applied Mathematical Sciences}, 139.
\newblock Springer-Verlag, New York, 1999.


\bibitem{Taylor}
G. I. Taylor.
\newblock The instability of liquid surfaces when accelerated in a direction perpendicular to their planes I.
\newblock {\em  Proc. Roy. Soc. London} A 201 (1950) 192-196.


\bibitem{TaylorBook}
M. Taylor.
\newblock Tools for PDE.
Pseudodifferential operators, paradifferential operators, and layer potentials. Mathematical Surveys and Monographs, 81.
\newblock {\em American Mathematical Society, Providence, RI}, 2000. x+257 pp.

\bibitem{TotzWuNLS}%[Wu12]
N. Totz and S. Wu.
\newblock A rigorous justification of the modulation approximation to the 2D full water wave problem.
\newblock {\em Comm. Math. Phys.} 310 (2012), no. 3, 817-883.

%\bibitem{Totz}
%N. Totz.
%\newblock modulation in 3D
%\newblock {\em }


%\bibitem{VBDias}
%J.-M. Vanden-Broeck and F. Dias.
%\newblock Gravity-capillary solitary waves in water of infinite depth and related free-surface flows.
%\newblock {\em J. Fluid Mech.} 240 (1992), 549-557.


\bibitem{Varvaruca}
E. Varvaruca.
\newblock Some geometric and analytic properties of solutions of Bernoulli free-boundary problems.
\newblock {\em Interfaces Free Bound.} 9 (2007), 367-381.


\bibitem{Wa1}
X. Wang.
\newblock Global infinite energy solutions for the 2D gravity water waves system.
\newblock Preprint {\em arXiv: 1502.00687}.


%%% REF % Checks

\bibitem{Wa2}
X. Wang.
\newblock On 3D water waves system above a flat bottom.
\newblock Preprint {\em arXiv:1508.06223}.

\bibitem{Wa3}
X. Wang.
\newblock Global solution for the 3D gravity water waves system above a flat bottom.
\newblock Preprint {\em arXiv:1508.06227}.

\bibitem{Wa4}
X. Wang.
\newblock Global regularity for the 3D finite depth capillary water waves.
\newblock Preprint {\em arXiv:1611.05472}.



%%% REF
%%% Numerics Wilkening

\bibitem{Wilkening1}
J. Wilkening.
\newblock Breakdown of Self-Similarity at the Crests of Large-Amplitude Standing Water Waves.
\newblock {\em Phys. Rev. Lett.} (2011) 107, 184501.

\bibitem{Wilkening2} % Put in
J. Wilkening and V. Vasal.
\newblock Comparison of five methods of computing the Dirichlet-Neumann operator for
the water wave problem.
\newblock {\em Nonlinear wave equations: analytic and computational techniques}
 175-210, Contemp. Math., 635, Amer. Math. Soc., Providence, RI, 2015.

%%%%%%
%%%%%%


\bibitem{Wu1}
S. Wu.
\newblock Well-posedness in Sobolev spaces of the full water wave problem in 2-D.
\newblock {\em Invent. Math.} 130 (1997), 39-72.


\bibitem{Wu2}
S. Wu.
\newblock Well-posedness in Sobolev spaces of the full water wave problem in 3-D.
\newblock {\em J. Amer. Math. Soc.} 12 (1999), 445-495.

\bibitem{WuVS}
S. Wu.
\newblock Mathematical analysis of vortex sheets.
\newblock {\em Comm. Pure Appl. Math.} 59 (2006), no. 8, 1065-1206.

\bibitem{WuAG}
S. Wu.
\newblock Almost global wellposedness of the 2-D full water wave problem.
\newblock {\em Invent. Math.} 177 (2009), 45-135.


\bibitem{Wu3DWW}
S. Wu.
\newblock Global wellposedness of the 3-D full water wave problem.
\newblock {\em Invent. Math.} 184 (2011), 125-220.


\bibitem{WuAC}
S. Wu.
\newblock A blow-up criteria and the existence of 2d gravity water waves with angled crests.
\newblock Preprint {\em arXiv:1502.05342}.

\bibitem{Wusurvey}
S. Wu.
\newblock Wellposedness and singularities of the water wave equations.
\newblock {\em Lectures on the theory of water waves}, 171-202,
\newblock {\em London Math. Soc. Lecture Note Ser.}, 426, Cambridge Univ. Press, Cambridge, 2016.



\bibitem{Yosi}
H. Yosihara.
\newblock Gravity waves on the free surface of an incompressible perfect fluid of finite depth.
\newblock {\em Publ. Res. Inst. Math. Sci.} 18 (1982), 49-96.


\bibitem{Zak0} %Hamiltonian formulation ??
V. E. Zakharov.
\newblock Stability of periodic waves of finite amplitude on the surface of a deep fluid.
\newblock {\em Zhurnal Prikladnoi Mekhaniki i Teckhnicheskoi Fiziki} 9 (1968), no.2,  86-94.
\newblock {\em J. Appl. Mech. Tech. Phys.}, 9, 1990-1994.


%\bibitem{Zak1}
%V. E. Zakharov
%\newblock
%\newblock

%%% REF %Check content

\bibitem{Zak2}
V. E. Zakharov.
\newblock Weakly nonlinear waves on the surface of an ideal finite depth fluid.
\newblock {\em Nonlinear waves and weak turbulence}, 167-197, Amer. Math. Soc. Transl. Ser. 2, 182, Amer. Math. Soc., Providence, RI, 1998.


\bibitem{ZakDya}
V. E. Zakharov and A. I. Dyachenko.
\newblock Is free-surface hydrodynamics an integrable system?
\newblock {\em Phys. Lett. A} 190 (1994), no. 2, 144-148.

\bibitem{ZakDya2}
V. E. Zakharov and A. I. Dyachenko.
\newblock High-Jacobian approximation in the free surface dynamics of an ideal fluid.
\newblock Nonlinear phenomena in ocean dynamics (Los Alamos, NM, 1995).
\newblock {\em Phys. D} 98 (1996), no. 2-4, 652-664.

\bibitem{ZhangZhang}
P. Zhang and Z. Zhang.
\newblock On the free boundary problem of three-dimensional incompressible Euler equations.
\newblock {\em Comm. Pure Appl. Math.} 61 (2008), no. 7, 877-940.


\end{thebibliography}
\end{document}